\documentclass[12pt]{amsart}
\usepackage{amscd}
\usepackage[dvips]{graphicx}
\usepackage{amscd}
\usepackage{latexsym}
\usepackage{amssymb}

\usepackage{color}

\newtheorem{theorem}{Theorem}[section]
\newtheorem{proposition}[theorem]{Proposition}
\newtheorem{prop}[theorem]{Proposition}
\newtheorem{lemma}[theorem]{Lemma}
\newtheorem{corollary}[theorem]{Corollary}
\newtheorem{defprop}[theorem]{Definition and Proposition}
\newtheorem{remark}[theorem]{Remark}
\newtheorem{definition}[theorem]{Definition}
\newtheorem{notation}[theorem]{Notation}

\newtheorem{example}[theorem]{Example}

\newcommand{\prf}{{\it Proof.} }

\newcommand{\qedd}{\hspace*{\fill}\hbox{$\Box$}}

\newcommand{\Zah}{{\Bbb Z}}

\newcommand{\mat}[2]{\mbox{${\rm M}_{#1}(#2)$}}

\newcommand{\tr}{\mbox{${\rm tr}$}}

\def\Spec{\mathop{\rm Spec}}

\begin{document}
\begin{center}{\large\bf The moduli of representations 
with Borel mold}\end{center}
\begin{center}{Kazunori NAKAMOTO\footnote{The author  
was partially supported by Grant-in-Aid in Scientific Research (C) (No. 23540044) 
from JSPS.}
 \\ {Center for Medical Education and Sciences, Faculty of Medicine}   
\\ University of Yamanashi \\ nakamoto@yamanashi.ac.jp}
\end{center}

\bigskip 

\begin{quote}\scriptsize 
{\bf Abstract.} 
The author constructs the moduli of representations whose images generate the subalgebra of upper triangular 
matrices (up to inner automorphisms) of the full matrix ring for any groups and any monoids. 
\end{quote} 

\begin{quote}\scriptsize 
\textup{2010} {\it Mathematics Subject Classification}. 
Primary 14D22; Secondary 14D20, 20M30, 20C99, 16G99.
\end{quote}

\begin{quote}\scriptsize 
\end{quote}

\bigskip 

%
%

The moduli of representations is very useful for studying 
representations and for describing the moduli spaces of various 
geometric objects.
There are several styles for constructing 
the moduli of representations.
One way is taking the quotient of the whole of 
the representation variety by 
${\rm PGL}$ (for example, see \cite{Simpson2}, \cite{kj-Saito93}, \cite{Nkmt00}, and so on).  
This method has a weak point: 
in \cite{Nkmt00} we constructed the coarse moduli scheme 
of equivalence classes of absolutely irreducible representations
as the universal geometric quotient of the open subset
consisting of stable points in the representation variety.  
However, on the complement of 
the absolutely irreducible representations, two representations 
which have the same composition factors become one 
point in the moduli of equivalence classes of representations. 
When two distinct representations have the same invariants, 
we can not separate them in the moduli.
If we want to separate two distinct representations,
we must choose another style.

In this paper we propose another style for constructing 
the moduli of representations in the non-absolutely irreducible 
case.
We introduce the notion of ``mold". 
A mold is, so to say, a subalgebra of the full matrix ring. 
We say that two representations have the same mold 
if their images generate the same type of subalgebras of the full matrix ring.
By using the notion of mold, 
we collect representations which have the same mold, 
and we construct the moduli of representations with a 
fixed mold. 
As an example of molds, we consider a Borel mold, 
that is, the subalgebra of 
upper triangular matrices (up to inner automorphisms). 
The main purpose of this article is the
construction of the moduli of (equivalence classes of) 
representations with Borel mold. 
This article is the first to develop ``mold program", 
that is, the construction of moduli schemes of representations 
from a viewpoint of mold.  
In another papers,  
we will construct several moduli schemes 
of representations with various molds. 

In \cite{Nkmt-deg2} we will deal with several molds 
of degree $2$.  
In the degree $2$ case, each molds over a field 
can be classified into $6$ types. 
The moduli of representation with Borel mold is one of 
$6$ types of the moduli of representations. 
In \cite{topos}, we have calculated the 
cohomology ring of the moduli of representations with 
Borel mold for free monoids over ${\mathbb C}$. 
Representations with Borel mold for free monoids are  
a special type of indecomposable modules over free algebras. 
The moduli of representations with 
Borel mold is the easiest geometric object to be investigated 
among the moduli spaces of indecomposable modules over free algebras. 
In \cite{rational}, we have investigated rational homotopy types of the moduli of representations with Borel mold 
for free monoids over ${\mathbb C}$ and related varieties.  
In \cite{notestopos}, we have studied the representation varieties with Borel mold 
for the case that the number of generators of the free monoid is small in comparison with the degree 
of representations. 

It is interesting and important to investigate 
the moduli of representations with Borel mold for 
several groups and monoids. 
Let us give an interesting example here. 
Let ${\rm Rep}_n(\Gamma)_B$ be the representation 
variety with Borel mold, that is, 
the subscheme consisting of representations with Borel 
mold in the representation variety.
Let us consider the universal representation with Borel mold on 
${\rm Rep}_n(\Gamma)_B$. 
The universal representation induces 
the action on the trivial vector bundle 
${\mathcal O}_{{\rm Rep}_n(\Gamma)_B}^{\oplus n}$ on ${\rm Rep}_n(\Gamma)_B$.
Then there exists a unique filtration of $\Gamma$-invariant 
subbundles $0={\mathcal E}_0 \subset {\mathcal E}_1 
\subset \cdots \subset {\mathcal E}_n = {\mathcal O}_{{\rm Rep}_n(\Gamma)_B}^{\oplus n}$ 
with ${\rm rk} \; {\mathcal E}_i =i$. 
When $n=2$ and $\Gamma$ is a free monoid (or a free group) 
of rank $2$, 
the universal sub-line bundle ${\mathcal E}_1$ is not
trivial on ${\rm Rep}_2(\Gamma)_B$, however 
${\mathcal E}_1^{\otimes 2}$ is trivial. 
From this fact, we see that 
each $2$-dimensional representation with Borel mold 
of a group generated by two elements 
on ${\rm Spec} R$ with $({\rm Pic} \:({\rm Spec} R))_2 = 0$ 
can be normalized into a representation in 
upper triangular matrices (Corollary \ref{cor-grt}). 
Here we denote by $({\rm Pic} \:({\rm Spec} R))_2$ the $2$-torsion part of the Picard group 
${\rm Pic} \:({\rm Spec} R)$. 
This fact shows one of geometric aspects of 
representations on schemes. 

By ``global representation theory" we understand 
a theory of representations on schemes, while 
by ``local representation theory" 
we understand a theory of representations on fields or local 
rings.  
"Global representation theory" has several geometric 
aspects like Corollary \ref{cor-grt}.    
The authors hopes that this article contributes to 
development of ``global representation theory". 

\bigskip

The author thanks Professor Akira Ishii for his useful advices on algebraic geometry. 
The author also would like to express thanks to Professor Takeshi Torii
for his technical advices on algebraic topology. 
The author thanks the referee for pointing out 
several mistakes and suggesting several ideas. 


%
%
%

\section{mold}
In this section, we introduce the notion of {\it mold}.
This notion is used for classification of representations and 
for constructing the moduli of representations.
We collect representations which have the 
same mold, and we attach a canonical scheme structure on the 
collection.
From a viewpoint of invariant theory, it is natural to 
classify representations with respect to mold.

\begin{definition}\rm
  Let $X$ be a scheme. A subsheaf of ${\mathcal O}_X$-algebras 
${\mathcal A} \subseteq {\rm M}_n({\mathcal O}_X)$ is said to be
a {\it mold} of degree $n$ on $X$ if 
${\mathcal A}$ and ${\rm M}_n({\mathcal O}_X)/{\mathcal A}$
are locally free sheaves on $X$.
We denote by ${\rm rk}{\mathcal A}$ the rank of ${\mathcal A}$ 
as a locally free sheaf on $X$. 
For a commutative ring $R$, we say that an $R$-subalgebra 
$A \subseteq {\rm M}_n(R)$ is a {\it mold} of degree $n$ 
over $R$ if 
$A$ is a mold of degree $n$ on $\Spec{R}$.
\end{definition}

We introduce the moduli of molds, that is, 
the moduli of subalgebras of the full matrix ring as follows:

\begin{defprop}
  The following contravariant functor is representable by 
a closed subscheme of the Grassmann scheme.
\[
\begin{array}{ccccl}
  {\mathcal Mold}_{n, d} & : & ({\bf Sch})^{op} & \to & ({\bf Sets}) \\
  & & X & \mapsto & \{ {\mathcal A} \mid 
\mbox{ a mold of deg } n \mbox{ on } 
X \mbox{ with } {\rm rk}{\mathcal A} = d \}.
\end{array}
\]
We denote by ${\rm Mold}_{n, d}$ the scheme representing the functor 
${\mathcal Mold}_{n, d}$. 
\end{defprop}

\prf Let ${\rm Grass}(d, {\rm M}_n)$ be the Grassmann scheme of 
rank $d$ subbundles of ${\rm M}_n$. 
The condition that a subbundle ${\mathcal A} \subset {\rm M}_n$ is closed under the multiplication 
of ${\rm M}_n$ and that ${\mathcal A}$ has the identity matrix is a closed condition. 
Hence the functor ${\mathcal Mold}_{n, d}$ is representable by the closed subscheme of 
${\rm Grass}(d, {\rm M}_n)$ defined by the condition above. 
\qed

\bigskip

We give some examples of the moduli of molds. 

\begin{example}\rm
 In the case $n = 2$, we have
\begin{eqnarray}
  {\rm Mold}_{2, 1}&  = & \Spec{{\Bbb Z}}, \label{ex1-1} \\
  {\rm Mold}_{2, 2}&  = & {\Bbb P}^2_{\Bbb Z}, \label{ex1-2} \\
  {\rm Mold}_{2, 3}&  = & {\Bbb P}^1_{\Bbb Z}, \label{ex1-3} \\
  {\rm Mold}_{2, 4}&  = & \Spec{\Bbb Z}.  \label{ex1-4}
\end{eqnarray}
Indeed, (\ref{ex1-1}) and (\ref{ex1-4}) are obvious.
To see (\ref{ex1-2}), note that giving an $R$-valued point of 
${\rm Mold}_{2, 2}$ is equivalent to giving 
a rank $1$ projective submodule of 
${\rm M}_2(R)/R\cdot I_2$ for each commutative ring $R$. 
Hence we have ${\rm Mold}_{2, 2} = {\Bbb P}^2_{\Bbb Z}$.
Later (\ref{ex1-3}) will be proved in Corollary 
\ref{cormaxp}.
\end{example}

We introduce an equivalence relation among  
molds as follows.

\begin{definition}\rm
  Let ${\mathcal A}$ and ${\mathcal B}$ be molds of degree 
$n$ on a scheme $X$. 
We say that ${\mathcal A}$ and ${\mathcal B}$ are 
{\it locally equivalent} if for each $x \in X$ 
there exist an neighborhood $U$ of $x$ and 
$P_{U} \in {\rm GL}_n({\mathcal O}_U)$ such that 
$P_{U}^{-1}({\mathcal A}\mid_U)P_{U} = {\mathcal B}\mid_U \subseteq 
{\rm M}_n({\mathcal O}_U)$.
\end{definition}

We define the following typical molds, a Borel mold and a 
parabolic mold.

\begin{definition}\rm
  We define the mold  ${\mathcal B}_n$ of degree $n$ on $\Spec{\Bbb Z}$
by
\[
{\mathcal B}_n := \{ (b_{ij}) \in {\rm M}_n({\Bbb Z}) \mid 
b_{ij} =0 \mbox{ for } i > j \}.
\]
Let ${\mathcal A}$ be a mold of degree $n$ on a scheme $X$.
We say that ${\mathcal A}$ is a {\it Borel mold} of degree $n$ 
if ${\mathcal A}$ and ${\mathcal B}_n\otimes_{\Bbb Z}{\mathcal O}_X$ 
are locally equivalent. 
\end{definition}

\begin{definition}\rm
Let $n_1, n_2, \ldots, n_r$ be positive integers with 
$\sum n_i = n$. 
We define the mold ${\mathcal P}_{n_1, n_2, \ldots, n_r}$ of 
degree $n$ on $\Spec{\Bbb Z}$ by
\[
{\mathcal P}_{n_1, n_2, \ldots, n_r} := 
\left\{ 
\begin{array}{c|c} 
    (b_{ij}) \in {\rm M}_n({\Bbb Z}) &  
  \begin{array}{c}
    b_{ij} =0 \mbox{ if } \sum_{k=1}^{s} n_k < i \le \sum_{k=1}^{s+1} n_k  \\
    \mbox{ and } j \le \sum_{k=1}^{s}n_k 
  \end{array} \\ 
\end{array}
\right\}.
\]  
Let ${\mathcal A}$ be a mold of degree $n$ on a scheme $X$.
We say that ${\mathcal A}$ is a {\it parabolic mold} of type 
$(n_1, n_2, \ldots, n_r)$
if ${\mathcal A}$ and ${\mathcal P}_{n_1, n_2, \ldots, n_r}
\otimes_{\Bbb Z}{\mathcal O}_X$ 
are locally equivalent. 
\end{definition}

\bigskip

Let us discuss the structure of 
the moduli of molds ${\mathcal Mold}_{n, d}$.
The following case is easy.

\begin{proposition}
For a positive integer $n$, we have 
\begin{eqnarray*}
  {\rm Mold}_{n, n^2} & = & \Spec{{\Bbb Z}}, \\
  {\rm Mold}_{n, d} & = & \emptyset \hspace*{10ex}
\mbox{ if } \quad n^2-n+1 < d < n^2.  
\end{eqnarray*}
\end{proposition}

\prf
Since there is no rank $n^2$ mold of degree $n$
except ${\rm M}_n$, we have ${\rm Mold}_{n, n^2}  =  \Spec{{\Bbb Z}}$.
Suppose that $n^2-n+1 < d < n^2$ and that 
$A \subseteq {\rm M}_n(k)$ is a rank $d$ mold over an algebraically 
closed field $k$.
Then $A$ has a non-trivial invariant subspace of $k^n$, 
and hence $A$ has at most dimension $n^2-n+1$.
This is a contradiction. Because there exists no 
geometric point of ${\rm Mold}_{n, d}$ if $n^2-n+1 < d < n^2$, 
we obtain ${\rm Mold}_{n, d}  =  \emptyset$. 
\qed

\bigskip

Let $n_1, n_2, \ldots, n_r$ be positive integers.
Put $n :=\sum_{1 \le i \le r} n_i$ and $d := \sum_{1 \le i \le j \le r} n_in_j$.
We show that the moduli of molds ${\rm Mold}_{n, d}$ contains 
an open and closed subscheme 
corresponding to the parabolic molds of type 
$( n_1, n_2, \ldots, n_r )$.
This subscheme is isomorphic to a flag scheme, and hence it it 
smooth over ${\Bbb Z}$.
To prove this statement, we make several preparations.  

\begin{notation}\rm
Let $n_1, n_2, \ldots, n_r$ be positive integers.

Put $n := \sum n_i$.
We define  
the closed subgroup scheme $P_{n_1, n_2, \ldots, n_r}$ of 
${\rm PGL}_n$ by 
\[
P_{n_1, n_2, \ldots, n_r} := 
\left\{ 
\begin{array}{c|c}
(b_{ij}) \in {\rm PGL}_n & 
\begin{array}{cc}
b_{ij}=0 & 
\mbox{ if }\sum_{\ell = 1}^{s} n_{\ell} < i \le 
\sum_{\ell =1}^{s+1} n_{\ell} \\  
& \mbox{ and } j \le 
\sum_{\ell = 1}^{s}n_{\ell}
\end{array} \\
\end{array}
\right\}.
\]
We denote by ${\rm Flag}_{n_1, n_2, \ldots, n_r}$ the 
flag scheme ${\rm PGL}_n/P_{n_1, n_2, \ldots, n_r}$.
\end{notation}

\begin{lemma}\label{lemma-mold-sub}
Let $R$ be a local ring.
Let us consider the canonical action of 
the parabolic mold ${\mathcal P}_{n_1, n_2, \ldots, n_r}\otimes_{\Zah}R$
on $R^n$ with $n = \sum n_i$. 
Then for each $1 \le s \le r$
there exists a unique 
rank $n_1 + n_2 + \cdots + n_s$ 
subbundle of $R^n$ which is 
invariant under the parabolic mold.  
{\rm (}By a {\it subbundle} $M$ of $R^n$ 
we understand an $R$-projective submodule 
$M$ of $R^n$ 
such that $R^n/M$ is also projective.{\rm )} 
\end{lemma}

{\it Proof.}
It is obvious that there exists an invariant 
rank $n_1+\cdots +n_s$ subbundle of $R^n$. 
For proving the uniqueness, we only have to show that 
the Borel mold ${\mathcal B}_n\otimes_{\Zah}R$ 
has a unique 
invariant rank $i$ subbundle of $R^n$ for each $1 \le i \le n$. 
Let $\{ e_1, e_2, \ldots, e_n \}$ be the canonical basis of $R^n$.
Suppose that  $M \subseteq R^n$ is an invariant rank $i$ subbundle.
Then we show that $M = Re_1 + Re_2 + \cdots Re_i$. 
If $v = \sum a_j e_j \in M$ with $a_j \in R^{\times}$ and $j > i$, 
then $\{ E_{1j}v, E_{2j}v, \ldots, E_{jj}v \}$ spans a 
rank $j$ subbundle of $M$. This is a contradiction.  
Hence if $v = \sum a_j e_j \in M$, then 
$a_j \in m$ for each 
$j > i$, where $m$ is the unique maximal ideal of $R$.

Let us define the projection $p : R^n \to R^{i}$ by 
$\sum a_j e_j \mapsto a_1 e_1 + \cdots + a_{i} e_{i}$.  
Since $M$ is a rank $i$ subbundle of $R^n$, 
$p\!\mid_M : M \to R^i$ is an isomorphism. 
For $1 \le j \le i$, we put $x_j := (p\!\mid_M)^{-1}(e_j)$. 
We can write $x_j = e_j + v_j$ with $v_j \in m e_{i+1} + \cdots + 
m e_{n}$. 
Then we have $E_{jj}x_j = e_{j} = (p\!\mid_M)^{-1}(e_{j}) 
= x_j = e_j + v_j$, which implies that $v_j = 0$.  
Therefore $M = Re_1 + Re_2 + \cdots Re_i$. 
\qed

\begin{corollary}\label{lemmanorma2}
Let $R$ be a local ring.
Then the set 
${\mathcal N}({\mathcal P}_{n_1, \ldots, n_r}\otimes_{\Zah}R) := \{ 
Q \in {\rm PGL}_n(R) \mid Q\cdot 
({\mathcal P}_{n_1, \ldots, n_r}\otimes_{\Zah}R)
\cdot Q^{-1}
\subseteq {\mathcal P}_{n_1, \ldots, n_r}\otimes_{\Zah}R \}$ is 
equal to $P_{n_1, \ldots, n_r}(R)$.
\end{corollary}

{\it Proof.}
Let $Q \in {\mathcal N}({\mathcal P}_{n_1, \ldots, n_r}\otimes_{\Zah}R)$.
Take a representative of $Q$ in ${\rm GL}_n(R)$, 
say it $Q$, too. 
Since the parabolic mold
${\mathcal P}_{n_1, \ldots, n_r}\otimes_{\Zah}R = Q\cdot 
({\mathcal P}_{n_1, \ldots, n_r}\otimes_{\Zah}R)
\cdot Q^{-1}$
has a unique invariant rank $n_1 + \cdots + n_s$ subbundles of $R^n$, 
$Q$ also leaves such a subbundle invariant. 
Hence $Q \in P_{n_1, \ldots, n_r}(R)$. 
\qed

\begin{prop}\label{prop:actionfree} 
Let $R$ be a local ring. 
For $Q \in P_{n_1, \ldots, n_r}(R)$, we define the 
algebra homomorphism ${\rm Ad}(Q) : 
{\mathcal P}_{n_1, \ldots, n_r}\otimes_{\Zah}R \to 
{\mathcal P}_{n_1, \ldots, n_r}\otimes_{\Zah}R$ by 
$Ad(Q)(X) = QXQ^{-1}$. 
If ${\rm Ad}(Q) = id$, then $Q = I_n$ in 
$P_{n_1, \ldots, n_r}(R)$. 
\end{prop}

\prf
Let $Q=(q_{ij}) \in P_{n_1, \ldots, n_r}(R)$. 
Suppose that ${\rm Ad}(Q)=id$. 
Let us consider 
the block $Q_k := (q_{ij})_{\sum_{m=1}^{k-1} n_m < i, j \le \sum_{m=1}^{k} n_m}$  
for $1 \le k \le r$. 
By the hypothesis, ${\rm Ad}(Q_k) : {\rm M}_{n_k}(R) \to 
{\rm M}_{n_k}(R)$ is the identity for $1 \le k  \le r$. 
Then we see that $Q_k$ is a scalar matrix. 
For $i < j$, considering the $(i, j)$-entry of $QE_{ij}Q^{-1} = E_{ij}$, we have $q_{ii}/q_{jj} = 1$. 
Hence $q_{11} = q_{22} = \cdots = q_{nn}$. 
For $i < j$, considering the $(i, j)$-entry of $QE_{jj}Q^{-1} = E_{jj}$, we have 
$q_{ij}/q_{jj} = 0$. Hence we obtain $q_{ij}=0$. Therefore $Q = I_n$ in  
$P_{n_1, \ldots, n_r}(R)$. 
\qed

\bigskip 

We construct a closed subscheme of the moduli of molds in 
the next proposition.

\begin{proposition}\label{propclosed}
Let $n_1, n_2, \ldots, n_r$ be positive integers.
Put $n := \sum n_i$ and $d := \sum_{1 \le i \le j \le r} n_in_j$.  
We define $\phi : {\rm PGL}_n \to {\rm Mold}_{n, d}$ by 
$Q \mapsto Q({\mathcal P}_{n_1, n_2, \ldots, n_r}\otimes_{\Zah}{\mathcal O}_X)Q^{-1}$ for 
a $X$-valued point $Q$ of ${\rm PGL}_n$ with a scheme $X$.
Then the morphism $\phi$ induces the closed immersion 
${\rm Flag}_{n_1, n_2, \ldots, n_r} \to {\rm Mold}_{n, d}$.
As a closed subscheme ${\rm Flag}_{n_1, n_2, \ldots, n_r}$ corresponds to 
the parabolic molds of type $(n_1, n_2, \ldots, n_r)$.
\end{proposition}

{\it Proof.}
The morphism $\phi$ induces $\overline{\phi} : 
{\rm PGL}_n/P_{n_1, \ldots, n_r} = {\rm Flag}_{n_1, \ldots, n_r} 
\to {\rm Mold}_{n, d}$.  
We claim that $\overline{\phi}$ is a closed immersion.
First we show that $\overline{\phi}$ is a monomorphism.
Let $X$ be a scheme.
Let $P$ and $Q$ be $X$-valued points of ${\rm PGL}_n$.
Suppose that $\phi(P)$ and $\phi(Q)$ are same molds on $X$.
Since  
$P({\mathcal P}_{n_1, \ldots, n_r}\otimes{\mathcal O}_{X, x})P^{-1} 
= Q({\mathcal P}_{n_1, \ldots, n_r}\otimes{\mathcal O}_{X, x})Q^{-1}$ 
for each $x \in X$, 
$P^{-1}Q$ is contained in 
the normalizer 
$N({\mathcal P}_{n_1, \ldots, n_r}\otimes{\mathcal O}_{X, x})$ 
of ${\mathcal P}_{n_1, \ldots, n_r}\otimes{\mathcal O}_{X, x}$.
From Corollary \ref{lemmanorma2} we see that  
$P^{-1}Q \in {\rm P}_{n_1, \ldots, n_r}$ at each $x$ 
and hence that $P = Q$ in ${\rm PGL}_n/{\rm P}_{n_1, \ldots, n_r}$.
Therefore $\overline{\phi}$ is a monomorphism.
Next the morphism $\overline{\phi}$ is proper, 
since the scheme ${\rm PGL}_n/{\rm P}_{n_1, \ldots, n_r}$ is proper 
over ${\Bbb Z}$. 
Thus we have proved that $\overline{\phi}$ is a closed immersion.
\qed

\bigskip

The closed subscheme constructed above is also open in 
the moduli of molds. 
For proving this, we introduce the following propositions. 

\begin{lemma}\label{lemmanormalizer}
Let $R$ be a local ring and let $A \subseteq {\rm M}_n(R)$ 
be a parabolic mold over $R$.
Then the normalizer $N(A) :=
\{ X \in {\rm M}_n(R) \mid [X, Y] \in  A \mbox{  for each }
Y \in A \}$ is equal to $A$.
Here we define $[X, Y] := XY-YX$.
\end{lemma}

{\it Proof.}
Since $R$ is a local ring, 
by changing $A$ to $PA P^{-1}$ with a suitable matrix 
$P \in {\rm GL}_n(R)$ we may assume that 
\begin{eqnarray}
\hspace*{5ex}  A = \left\{ 
\begin{array}{c|c}
(b_{ij}) \in {\rm M}_n(R) & 
\begin{array}{cc}
b_{ij}=0 & \mbox{ if }\sum_{\ell =1}^{s} n_{\ell} < i \le 
\sum_{\ell =1}^{s+1} n_{\ell} \\ 
& \mbox{ and } j \le 
\sum_{\ell = 1}^{s}n_{\ell}   
\end{array} \\
\end{array} 
\right\}. \label{eqnnorma}
\end{eqnarray} 
It is clear that $A \subseteq N(A)$.
Suppose that $X = \sum a_{ij}E_{ij} \in {\rm M}_n(R)\setminus A$.
There exists $E_{ij} \notin A$ with $a_{ij} \neq 0$. 
Note that $i > j$.
For $E_{jj} \in A$, we have   
$[ X, E_{jj} ] = a_{ij}E_{ij} + \cdots$.
Since the $(i, j)$-entry of $[ X, E_{ij} ]$
is not zero, 
$[ X, E_{jj} ] \notin A$. 
Hence $X \notin N(A)$.
Thus we have proved that $A = N(A)$.
\qed

\bigskip

For a mold $A \subseteq {\rm M}_n(k)$ with 
a field $k$, 
$A$ and ${\rm M}_n(k)/A$ are 
$A$-bimodules. 
We define 
\[ 
{\rm Der}_k(A, {\rm M}_n(k)/A) := 
\left\{ \begin{array}{c|l} 
 &  \mbox{ $\delta$ is $k$-linear and } \\ 
\delta : A \to {\rm M}_n(k)/A &  \delta(ab) = a \delta(b) + \delta(a) b  \\ 
 & \mbox{ for $a, b \in A$ } 
\end{array} 
\right\}.  
\] 

\begin{proposition}\label{propder}
Let $k$ be a field and let $A \subseteq {\rm M}_n(k)$ 
be a parabolic mold over $k$.
Then the linear map  
\[
\begin{array}{ccc}
{\rm M}_n(k)/A & \to & {\rm Der}_k(A, {\rm M}_n(k)/A) \\
X & \mapsto & [ X, - ]
\end{array}
\]
is bijective.
\end{proposition}

{\it Proof.}
We can easily check that the above map is well-defined.
The injectivity of the map follows from 
Lemma \ref{lemmanormalizer}.
For proving that the linear map is an isomorphism, 
we may assume 
(\ref{eqnnorma}).
Let $\delta \in {\rm Der}_k(A, {\rm M}_n(k)/A)$.  
If $E_{ij} \in A$, then 
we have 
\[
\delta(E_{ij}) = \delta(E_{ii}E_{ij}) = 
E_{ii}\delta(E_{ij}) + \delta(E_{ii})E_{ij} 
\]
and 
\[
\delta(E_{ij}) = \delta(E_{ij}E_{jj}) = 
E_{ij}\delta(E_{jj}) + \delta(E_{ij})E_{jj}.
\]
The first equality shows that 
the $(\ell, \ast)$-entries of $\delta(E_{ij})$ 
are determined by $\delta(E_{ii})$ for $\ell \neq i$, and 
the second equality shows that 
the $(\ast, \ell)$-entries of $\delta(E_{ij})$ are  
determined by $\delta(E_{jj})$ for $\ell \neq j$.
Hence $\delta(E_{ij})$ is 
determined by $\delta(E_{ii})$ and $\delta(E_{jj})$.  
The derivation $\delta$ is determined by 
$\{ \delta(E_{ii}) \mid 1 \le i \le n \}$. 

Since $\delta(E_{ii}) = \delta(E_{ii}E_{ii}) = 
E_{ii}\delta(E_{ii}) + \delta(E_{ii})E_{ii}$, 
the matrix $\delta(E_{ii})$ has zero entries except 
$(i, *)$-entries and $(*, i)$-entries. 
For $i \neq j$, we obtain 
\[
0 = \delta(E_{ii}E_{jj}) = E_{ii}\delta(E_{jj}) + \delta(E_{ii})E_{jj}.
\] 
If $j \le \sum_{\ell = 1}^{s} n_{\ell} < i$ for some $s$, then 
$\delta(E_{ii})_{ij} = - \delta(E_{jj})_{ij}$.
We see that $\delta$ is determined by the data 
$\{ \delta(E_{jj})_{ij} \mid j \le \sum_{\ell = 1}^{s} n_{\ell} < i 
\mbox{ for some } s \}$ and that $\dim {\rm Der}_k(A, {\rm M}_n(k)/A) 
\le \sum_{1 \le i < j \le r}n_in_j = \dim{\rm M}_n(k)/A$.
From the injectivity,  
we prove that the above linear map is an isomorphism. 
\qed

\bigskip
Now we prove that the closed subscheme ${\rm Flag}_{n_1, \ldots, n_r}$ 
is open in ${\rm Mold}_{n, d}$.

\begin{proposition}\label{prop:opensmooth} 
The morphism $\overline{\phi} : {\rm Flag}_{n_1, \ldots, n_r} \to 
{\rm Mold}_{n, d}$ in Proposition \ref{propclosed} is smooth.
In particular, ${\rm Flag}_{n_1, \ldots, n_r}$ is an open and 
closed subscheme of ${\rm Mold}_{n, d}$.
\end{proposition}

{\it Proof.}
Let $(R, m, k)$ be an artin local ring. Let $I$ be 
an ideal with $m\cdot I = 0$. Suppose that 
$A \subseteq {\rm M}_n(R)$ is a mold over $R$ such that 
$A\otimes_{R}R/I$  is a parabolic mold of type $(n_1, \ldots, 
n_r)$ over $R/I$.
For proving the statement, we only have to show that
$A$ is a parabolic mold over $R$. 
From the assumption, 
there exists $\overline{P} \in {\rm GL}_n(R/I)$ such that 
$\overline{P}(A\otimes_{R}R/I)\overline{P}^{\:-1} = 
{\mathcal P}_{n_1, \ldots, n_r}\otimes_{\Zah}R/I$. 
Take a matrix $P \in {\rm GL}_n(R)$ such that 
$P$ mod $I = \overline{P}$. 
By changing $A$ to $PAP^{-1}$, we may assume that 
$A$ is a mold over $R$ such that $A\otimes_{R}R/I = 
{\mathcal P}_{n_1, \ldots, n_r}\otimes_{\Zah}R/I$.
We denote $A\otimes_{R}R/I$ by $\overline{A}$.

Let us denote ${\mathcal P}_{n_1, \ldots, n_r} 
\otimes_{{\mathbb Z}} M$ by ${\mathcal P}_{n_1, \ldots, n_r}(M)$ 
for an $R$-module $M$.   
Let $q : {\rm M}_n(R) \to {\rm M}_n(R)/{\mathcal P}_{n_1, \ldots, n_r}(R)$ 
be the canonical projection. 
For matrix elements $E_{ij} \in \overline{A}$ 
choose their representatives $\tilde{E}_{ij} \in A$ 
in Lemma \ref{lemma:representative}. 
We define the $R/I$-linear map 
$\delta : \overline{A} \to {\rm M}_n(I)/{\mathcal P}_{n_1, \ldots, n_r}(I)$ by 
$\delta( \sum {a}_{ij} E_{ij} ) := 
q( \sum \tilde{a}_{ij} \tilde{E}_{ij} )$, where 
$\tilde{a}_{ij} \in R$ is a representative of $a_{ij} \in R/I$. 
Since $I^2 = 0$, we can easily check that 
$\delta$ is independent of choices of 
$\tilde{a}_{ij}$. 

Note that 
${\rm M}_n(I) = {\rm M}_n(R)\otimes_{R}I\otimes_{R}R/m 
= {\rm M}_n(k)\otimes_{k}I$.  
The map $\delta$ induces the $k$-linear map 
$\overline{\delta} : A\otimes_{R}k = \overline{A}\otimes_{R/I}k 
\to ({\rm M}_n(k)/{\mathcal P}_{n_1, \ldots, n_r}(k))\otimes_{k}I$. 
Let us show that 
$\overline{\delta}$ is a derivation. 
Denote by $E'_{ij} \in A$ the matrix elements in ${\rm M}_n(R)$. 
Set $\tilde{E}_{ij} = E'_{ij} + X_{ij}$, where  
$X_{ij} := \tilde{E}_{ij} - E'_{ij} \in {\rm M}_n(I)$. 
Then 
$\tilde{E}_{ij}\tilde{E}_{k\ell} = (E'_{ij} + X_{ij})(E'_{k\ell} + X_{k\ell}) = \delta_{jk} E'_{i\ell} 
+ E'_{ij}X_{k\ell} + X_{ij}E'_{k\ell}$.  
On the other hands, $\tilde{E}_{ij}\tilde{E}_{k\ell} = \delta_{jk} \tilde{E}_{i\ell} 
+\sum_{(s, t) \neq (i, \ell)} x_{st} \tilde{E}_{st}$ with some $x_{st} \in I$ 
because $\oplus_{(k, \ell) \in {\mathcal P}} R \tilde{E}_{k\ell} = A$. 
Hence $q(E'_{ij}X_{k\ell} + X_{ij}E'_{k\ell}) = q(\delta_{jk} E'_{i\ell} 
+ E'_{ij}X_{k\ell} + X_{ij}E'_{k\ell}) = q(\delta_{jk} \tilde{E}_{i\ell} 
+\sum_{(s, t) \neq (i, \ell)} x_{st} \tilde{E}_{st}) = \delta_{jk} q(\tilde{E}_{i\ell})$. 
Since 
$\overline{\delta}(E_{ij}E_{k\ell}) = \overline{\delta}(\delta_{jk}E_{i\ell}) 
= \delta_{jk} q(\tilde{E}_{i\ell})$
and 
$E_{ij}\overline{\delta}(E_{k\ell}) + \overline{\delta}(E_{ij})E_{k\ell} 
 =  E'_{ij}q(X_{k\ell}) + q(X_{ij})E'_{k\ell}   
 =  q(E'_{ij}X_{k\ell} + X_{ij}E'_{k\ell})$, we have 
$\overline{\delta}(E_{ij}E_{k\ell}) =  
E_{ij}\overline{\delta}(E_{k\ell}) + \overline{\delta}(E_{ij})E_{k\ell}$.   
Therefore $\overline{\delta}$ is a derivation. 

By Proposition \ref{propder} we have $Y \in {\rm M}_n(k)\otimes_{k}I = 
{\rm M}_n(I)$ such that $\overline{\delta} = [Y, -]$. 
Hence the map $A \stackrel{{\rm proj.}}{\to} \overline{A} \stackrel{\delta}{\to} 
{\rm M}_n(I)/{\mathcal P}_{n_1, \ldots, n_r}(I)$ 
is given by $X \mapsto [Y, X] \mod {\mathcal P}_{n_1, \ldots, n_r}(I)$.  
Putting $P:=(I_n - Y) \in {\rm GL}_n(R)$, we have 
$PXP^{-1} = (I_n - Y)X(I_n + Y) = X - [Y, X]$ for $X \in {\rm M}_n(R)$.
We see that $PAP^{-1} = {\mathcal P}_{n_1, \ldots, n_r}\otimes_{\Zah}R$
and hence that $A$ is a parabolic mold over $R$. 
This completes the proof.   
\qed

\bigskip

The following lemma has been used in Proposition \ref{prop:opensmooth}. 

\begin{lemma}\label{lemma:representative}
Let $(R, m, k)$ be an artin local ring. 
Let $I$ be an ideal of $R$ with $m\cdot I = 0$. 
For $X = (x_{ij}) \in {\rm M}_n(R)$, ${\rm Supp} X := 
\{ (i, j) \mid x_{ij} \neq 0 \}$. 
Set ${\mathcal P} := 
\{ (i, j) \mid \sum_{\ell = 1}^{s} n_{\ell} < i \le \sum_{\ell = 1}^{s+1} n_{\ell} 
\mbox{ and } 
\sum_{\ell = 1}^{s} n_{\ell} < j \mbox{ for some } s
\}$. 
Assume that $A \subseteq {\rm M}_n(R)$ is a mold over $R$ such that 
$A\otimes_{R}R/I = {\mathcal P}_{n_1, \ldots, n_r}\otimes_{\Zah}R/I$. 
Then there exist $\tilde{E}_{ij} \in A$ for $(i, j) \in {\mathcal P}$ satisfying 
the following properties: 
\begin{enumerate} 
\item  $\tilde{E}_{ij} \mod I$ coincides with the matrix element $E_{ij} \in {\rm M}_n(R/I)$ 
and the $(i, j)$-entry of $\tilde{E}_{ij}$ is $1$ for each $(i, j) \in {\mathcal P}$. 
\item ${\rm Supp} \tilde{E}_{ij} \subseteq \{ (i, j) \} \cup 
{\mathcal P}^{c}$ 
for each $(i, j) \in {\mathcal P}$.  
\item $\oplus_{(k, \ell) \in {\mathcal P}} R \tilde{E}_{k\ell} = A$. 
\end{enumerate} 
\end{lemma} 
 
\prf
Take $\tilde{E}_{ij} \in A$ such that 
$\tilde{E}_{ij} \mod I$ coincides with $E_{ij} \in {\rm M}_n(R/I)$ and   
the $(i, j)$-entry of $\tilde{E}_{ij}$ is $1$. 
By changing $\tilde{E}_{ij}$ into $\tilde{E}_{ij} - \sum_{{\mathcal P} \ni (k, \ell) \neq (i, j)} x_{k\ell} \tilde{E}_{k\ell}$ 
for some $x_{k\ell} \in I$, 
we may assume that ${\rm Supp} \tilde{E}_{ij} 
\subseteq \{ (i, j) \} \cup {\mathcal P}^{c}$ 
for each $(i, j) \in {\mathcal P}$.     

Let us show that $\oplus_{(k, \ell) \in {\mathcal P}} R \tilde{E}_{k\ell} = A$. 
Indeed, $\{ \tilde{E}_{k\ell} \}_{(k, \ell) \in {\mathcal P}}$ are 
linearly independent over $R$ and there exist the following exact sequences:  
\[ 
\begin{array}{ccccccccc}
0 & \to & N &  
\to & {\rm M}_n(R) & \to &  {\rm M}_n(R)/N  
& \to & 0 \\ 
&  & \downarrow & & \parallel & & \downarrow & \\
0 & \to & A &  
\to & {\rm M}_n(R) & \to &  {\rm M}_n(R)/A   
& \to & 0,  
\end{array} 
\]
where $N := \oplus_{(k, \ell) \in {\mathcal P}} R \tilde{E}_{k\ell}$. 
Since ${\rm M}_n(R)/N$ and ${\rm M}_n(R)/A$ are isomorphic to 
free $R$-modules of rank $\sharp {\mathcal P}^c$, 
${\rm M}_n(R)/N \cong {\rm M}_n(R)/A$.  
Hence $N = A$. 
\qed

\bigskip 

From the discussion above, we obtain the following theorem. 

\begin{theorem}
Let $n_1, n_2, \ldots, n_r$ be positive integers.
Put $n := \sum n_i$ and $d := \sum_{1 \le i \le j \le r} n_in_j$.
Then the moduli of molds  ${\rm Mold}_{n, d}$ contains the open 
and closed subscheme corresponding to the 
parabolic molds of type $(n_1, n_2, \ldots, n_r)$.
This subscheme is isomorphic to a flag scheme over ${\Bbb Z}$.
\end{theorem}

\bigskip 

The above theorem follows the next corollary. 

\begin{corollary}\label{cormaxp}
In the case $d = n^2 - n + 1$, we have 
\[ {\rm Mold}_{n, n^2-n+1} \cong 
\left\{
\begin{array}{cc}
{\Bbb P}_{\Zah}^{n-1} \coprod {\Bbb P}_{\Zah}^{n-1} & n > 2 \\
{\Bbb P}_{\Zah}^1 & n=2. 
\end{array}
\right.
\]  
\end{corollary}

\prf
Each geometric point in ${\rm Mold}_{n, n^2-n+1}(\Omega)$, that is, 
each mold of rank $n^2 - n +1$ over $\Omega$  has 
an invariant subspace of $\Omega^n$. 
This mold is a parabolic mold of type $(1, n-1)$ or $(n-1, 1)$. 
Hence the moduli ${\rm Mold}_{n, n^2-n+1}$ is covered by 
${\rm Flag}_{1, n-1}$ and ${\rm Flag}_{n-1, 1}$. 
\qed

\bigskip

From now on, we prepare some terminologies on representations.  
Using the notion of mold, we classify representations. 

\begin{definition}\rm\label{def-mold-rep}
Let $\Gamma$ be a group or a monoid. Let $X$ be a scheme.
By a {\it representation} of degree $n$ on $X$  for $\Gamma$
we understand a group homomorphism (resp. a monoid homomorphism) 
$\rho : \Gamma \to {\rm GL}_n(\Gamma(X, {\mathcal O}_X))$ 
(resp. 
$\rho : \Gamma \to {\rm M}_n(\Gamma(X, {\mathcal O}_X))$),
where $\Gamma(X, {\mathcal O}_X)$ is the ring of global 
sections on $X$.

For two representations $\rho, \rho'$ of degree $n$ for $\Gamma$ 
on a scheme $X$, we say that $\rho$ and $\rho'$ are {\it equivalent} 
 (or $\rho \sim \rho'$) 
if there exists 
a $\Gamma(X, {\mathcal O}_X)$-algebra isomorphism 
$\sigma : {\rm M}_n(\Gamma(X, {\mathcal O}_X)) 
\to {\rm M}_n(\Gamma(X, {\mathcal O}_X))$ such that 
$\sigma(\rho(\gamma)) = \rho'(\gamma)$ for each $\gamma \in \Gamma$. 
We also say that $\rho$ and $\rho'$ are {\it locally equivalent}
if there exists an open covering $X = \cup_{i \in I}U_i$ such that 
$\rho\!\mid_{U_i}$ and $\rho'\!\mid_{U_i}$ are equivalent for each 
$i \in I$. 
\end{definition}

\begin{definition}\rm
Let ${\mathcal A}$ be a mold of degree $n$ on $X$.
For a representation $\rho$ on $X$, we say that 
$\rho$ is a {\it representation with mold} ${\mathcal A}$
if the subsheaf ${\mathcal O}_X[\rho(\Gamma)]$ 
of ${\mathcal O}_X$-algebras of
${\rm M}_n({\mathcal O}_X)$ 
is locally equivalent to ${\mathcal A}$.
In particular, we say that $\rho$ is a 
{\it representation with Borel mold} if 
${\mathcal O}_X[\rho(\Gamma)]$ is a Borel mold.
We also say that $\rho$ is a 
{\it representation with parabolic mold of type 
$(n_1, n_2, \ldots, n_r)$} if 
${\mathcal O}_X[\rho(\Gamma)]$ is a parabolic mold of type 
$(n_1, n_2, \ldots, n_r)$.
\end{definition}

In \cite{Nkmt00} we proved the existence of the coarse moduli scheme 
of equivalence classes of absolutely irreducible representations.
Here we quote this result.

\begin{definition}\rm
  For a representation $\rho$ of degree $n$ for a 
group or monoid $\Gamma$ on a scheme $X$, we say that 
$\rho$ is {\it absolutely irreducible} if 
${\mathcal O}_X[\rho(\Gamma)] = {\rm M}_n({\mathcal O}_X)$.
This definition is equivalent to the one in \cite{Nkmt00}.
We abbreviate an absolutely irreducible representation to a.i.r.
\end{definition}

\begin{theorem}[\cite{Nkmt00}]
There exists a coarse moduli scheme separated 
over ${\Bbb Z}$ associated to the 
following moduli functor:
\[
\begin{array}{ccccl}
{\mathcal E}q{\mathcal AIR}_n(\Gamma) & : & ({\bf Sch}) & 
\to & ({\bf Sets}) \\
 & & X & \mapsto & \{ \rho : \mbox{ a.i.r. of degree $n$ for $\Gamma$ 
on } X \}/\sim. 
\end{array}
\]
In particular, if $\Gamma$ is a finitely generated group (or monoid), 
then 
the moduli is of finite type over ${\Bbb Z}$.
\end{theorem}

In the sequel, 
we only deal with representations with Borel mold. 
We will construct the moduli schemes of equivalence classes  
of representations with Borel mold.

\section{Construction of the moduli of representations with 
Borel mold}

In this section, we construct the moduli scheme of 
equivalence classes of representations with Borel mold.

\bigskip

Let us recall the representation variety (For details, see \cite{Nkmt00}). 
Let $\Gamma$ be a group or a monoid.
The following contravariant functor 
is representable by an affine scheme: 
\[
\begin{array}{ccccl}
	{\rm Rep}_n(\Gamma) & : & ({\bf Sch})^{op} & \to & ({\bf Sets}) \\
	 & & X & \mapsto & \{ \rho : \mbox{ rep. of deg $n$ for $\Gamma$ 
	 on $X$ } \}. 
\end{array}
\]
We call the affine scheme ${\rm Rep}_n(\Gamma)$ the 
representation variety of degree $n$ for $\Gamma$. 
The group scheme ${\rm PGL}_n$ over ${\Bbb Z}$ 
acts on ${\rm Rep}_n(\Gamma)$ 
by $\rho \mapsto P^{-1}\rho P$.
Each ${\rm PGL}_n$-orbit forms an equivalence class of 
representations. 

\bigskip

For a commutative ring $R$, we set 
${\mathcal B}_n(R) := \{ (a_{ij}) \in {\rm M}_n(R) 
\mid a_{ij} = 0 \mbox{ for each } i>j \}$, that is, 
${\mathcal B}_n(R)$ is the $R$-subalgebra of 
upper triangular matrices.
We define the closed subgroup scheme ${\rm B}_n$ of 
${\rm PGL}_n$ by 
${\rm B}_n := \{ (a_{ij}) \in {\rm PGL}_n \mid 
a_{ij}=0 \mbox{ for } i>j \}$. 
By a {\it ${\mathcal B}_n$-representation} of degree $n$ for $\Gamma$ 
on a scheme $X$, 
we understand a homomorphism 
$\rho : \Gamma \to {\mathcal B}_n(\Gamma(X, {\mathcal O}_X))$. 
It is easy to check that the subfunctor of ${\rm Rep}_n(\Gamma)$ 
consisting of ${\mathcal B}_n$-representations 
is represented by a closed subscheme of ${\rm Rep}_n(\Gamma)$. 
We denote it by ${\rm B}_n(\Gamma)$.  

For two ${\mathcal B}_n$-representations 
$\rho$ and $\rho'$ of degree $n$ for $\Gamma$
on a scheme $X$, we say that 
$\rho$ and $\rho'$ are
{\it ${\rm B}_n$-equivalent} to each other, if 
there exists a $X$-valued point $Q \in {\rm B}_n(X)$ 
such that $Q\rho Q^{-1}= \rho'$. 
The group scheme ${\rm B}_n$ acts on ${\rm B}_n(\Gamma)$ 
by $\rho  \mapsto Q\rho Q^{-1}$. 
Each ${\rm B}_n$-orbit forms a ${\rm B}_n$-equivalence class of
${\mathcal B}_n$-representations. 

By a {\it ${\mathcal B}_n$-representation with Borel mold} for 
$\Gamma$ on a scheme $X$, we understand 
a homomorphism $\rho : \Gamma \to {\mathcal B}_n(\Gamma(X, 
{\mathcal O}_X))$ which is a representation with 
Borel mold. Note that 
$\rho : \Gamma \to {\mathcal B}_n(\Gamma(X, 
{\mathcal O}_X))$ is a ${\mathcal B}_n$-representation 
with Borel mold if and only if 
$\rho(\Gamma)$ generates ${\mathcal B}_n(\Gamma(X, 
{\mathcal O}_X))$ as a $\Gamma(X, {\mathcal O}_X)$-algebra. 
We denote by ${\rm Rep}_n(\Gamma)_B$ the subfunctor of 
${\rm Rep}_n(\Gamma)$ consisting of 
representations with Borel mold. 
We also denote by ${\rm B}_n(\Gamma)_B$ the 
subfunctor of ${\rm B}_n(\Gamma)$ consisting of 
${\mathcal B}_n$-representations with Borel mold.
We see that 
${\rm Rep}_n(\Gamma)_B$ and ${\rm B}_n(\Gamma)_B$ 
are subschemes of 
${\rm Rep}_n(\Gamma)$ and ${\rm B}_n(\Gamma)$, 
respectively. 
Indeed, we can verify that ${\rm Rep}_n(\Gamma)$ 
has a (locally closed) subscheme ${\rm Rep}_n(\Gamma)_d$ of representations with 
mold of rank $d = n(n+1)/2$. 
Then the subscheme ${\rm Rep}_n(\Gamma)_B$ can be obtained by 
taking the pull-back of ${\rm Flag}_{1, 1, \ldots, 1}$ by the morphism 
${\rm Rep}_n(\Gamma)_d \to {\rm Mold}_{n, d}$. 
Similarly, we can verify that ${\rm B}_n(\Gamma)_B$ is an 
open subscheme of ${\rm B}_n(\Gamma)$.  
We call the scheme ${\rm Rep}_n(\Gamma)_B$ 
the {\it representation variety with Borel mold} of degree 
$n$ for $\Gamma$. 

\bigskip

Now let us define the contravariant 
functor ${\mathcal EqB}_n(\Gamma)$. 
By a {\it generalized 
representation with 
Borel mold} of degree $n$  
for $\Gamma$ on a scheme $X$, we understand 
pairs $\{ (U_i, \rho_i) \}_{i \in I}$ of 
an open set $U_i$ and a representation with Borel mold 
$\rho_i : \Gamma \to {\rm M}_n(\Gamma(U_i, {\mathcal O}_X))$ 
satisfying the following two conditions: 
\begin{enumerate}
\item $\cup_{i \in I} \;U_i =X$,  \label{cond-gen-1} 
\item for each $x \in U_i \cap U_j$, $\rho_i$ and $\rho_j$ are 
equivalent to each other on a neighbourhood of $x$. 
\label{cond-gen-2}
\end{enumerate}

\bigskip 

We say that generalized representations with Borel mold  
$\{ (U_i, \rho_i ) \}_{i\in I}$ and $\{ (V_j, \sigma_j) \}_{j \in J}$ 
are {\it equivalent} (or $\{ (U_i, \rho_i ) \}_{i\in I} \sim  
\{ (V_j, \sigma_j) \}_{j \in J}$) if $\{ (U_i, \rho_i ) \}_{i\in I} \cup 
\{ (V_j, \sigma_j) \}_{j \in J}$ is a generalized representation 
again. 

\bigskip

We introduce the contravariant functor 
${\mathcal EqB}_n(\Gamma)$:   
\[
\begin{array}{ccccl}
{\mathcal EqB}_n(\Gamma) & : & ({\bf Sch})^{op} & \to & ({\bf Sets}) \\
  & & X & \mapsto & 
  \left\{ 
  \begin{array}{c|l} 
  & \mbox{ gen. rep. with } \\ 
  \{  (U_i, \rho_i) \}_{i \in I} 
  & \mbox{ Borel mold of } \\ 
   & \mbox{ deg $n$ for $\Gamma$ on $X$ } 
   \end{array}
   \right\} \Bigg/\sim.
\end{array}
\]
In this section, we show that 
the functor ${\mathcal EqB}_n(\Gamma)$ is 
representable by a scheme over ${\Bbb Z}$.
For proving this, we prepare another functor 
${\mathcal EqBB}_n(\Gamma)$. 

\bigskip

By a {\it generalized ${\mathcal B}_n$-representation} with 
Borel mold of degree $n$  
for $\Gamma$ on a scheme $X$, we understand 
pairs $\{ (U_i, \rho_i) \}_{i \in I}$ of 
an open set $U_i$ and a ${\mathcal B}_n$-representation 
with Borel mold 
$\rho_i : \Gamma \to {\mathcal B}_n(\Gamma(U_i, {\mathcal O}_X))$ 
satisfying the above condition (\ref{cond-gen-1}) and the following: \\

(\ref{cond-gen-2})$^{\ast}$ 
for each $x \in U_i \cap U_j$, there exists $Q \in {\rm B}_n(V)$ 
 such that $Q^{-1}\rho_i Q = \rho_j$ on a neighbourhood $V$ of 
$x$. 

\bigskip 

We say that two generalized ${\mathcal B}_n$-representations 
$\{ (U_i, \rho_i ) \}_{i\in I}$ and $\{ (V_j, \sigma_j) \}_{j \in J}$ 
are {\it ${\rm B}_n$-equivalent} 
(or $\{ (U_i, \rho_i ) \}_{i\in I} \sim_{{\rm B}_n} 
\{ (V_j, \sigma_j) \}_{j \in J}$)
if $\{ (U_i, \rho_i ) \}_{i\in I} \cup 
\{ (V_j, \sigma_j) \}_{j \in J}$ is a generalized 
${\mathcal B}_n$-representation  
again.

\bigskip 

We define the contravariant functor 
${\mathcal EqBB}_n(\Gamma)$ by  
\[
\begin{array}{ccccl}
{\mathcal EqBB}_n(\Gamma) & : & ({\bf Sch})^{op} & \to & ({\bf Sets}) \\
  & & X & \mapsto & 
  \left\{ 
  \begin{array}{c|l} 
   & \mbox{ gen. ${\mathcal B}_n$-rep. with } \\ 
  \{ (U_i, \rho_i) \}_{i \in I} 
  & 
   \mbox{ Borel mold of } \\ 
   & \mbox{ deg $n$ for $\Gamma$ on $X$ } 
   \end{array}
   \right\} \Bigg/\sim_{{\rm B}_n}.
\end{array}
\]

\bigskip

We can check that  
there exists a canonical isomorphism 
${\mathcal EqBB}_n(\Gamma) \to 
{\mathcal EqB}_n(\Gamma)$. 
Hence the representability of ${\mathcal EqB}_n(\Gamma)$ 
is reduced to the one of ${\mathcal EqBB}_n(\Gamma)$. 

\bigskip 

The following lemma can be easily verified: 

\begin{lemma}
The functor ${\mathcal EqBB}_n(\Gamma)$ 
is a sheaf with respect to Zariski topology. 
\end{lemma}

\bigskip

By the above lemma, 
for proving that ${\mathcal EqBB}_n(\Gamma)$ is representable,
it suffices to show that 
it admits an open covering of  
affine schemes. 

\bigskip

Let us consider the index set ${\mathcal I}_n := \{ (i, j) \mid 
1 \le i \le j \le n \}$. 
We define the order on ${\mathcal I}_n$ by 
\[ 
(i, j) \le (i', j') \Leftrightarrow \left\{
\begin{array}{l} 
|i-j| < |i'-j'|  \\
\mbox{ or } \\
|i-j| = |i'-j'|  \mbox{ and } i \le i' .
\end{array}
\right. 
\]
We define $J_{i, j}(R) := \{ (a_{kl}) \in {\mathcal B}_n(R) 
\mid a_{kl}=0 \mbox{ for each }  (k, l) \le (i, j) \}$.

Let $\rho$ be a ${\mathcal B}_n$-representation with Borel mold 
of degree $n$ for $\Gamma$ 
on a scheme $X$. 
We say that $\rho$ satisfies the $(\ast)$-condition with 
respect to 
an ${\mathcal I}_n$-indexed 
subset $\{ \alpha_{ij} \}_{(i,j) \in {\mathcal I}_n}$ of 
$\Gamma$ 
if the set 
$\{ \rho(\alpha_{k, \ell}) \}_{(k, \ell) \le (i, j)}$ 
forms a basis of ${\mathcal B}_n(\Gamma(X, {\mathcal O}_X))
/J_{i, j}(\Gamma(X, {\mathcal O}_X))$ over 
$\Gamma(X, {\mathcal O}_X)$ 
for each $(i, j) \in {\mathcal I}_n$.
For any $P \in {\rm B}_n(X)$, 
$\rho$ satisfies the $(\ast)$-condition 
with respect to 
$\{ \alpha_{i, j} \}_{(i, j) \in {\mathcal I}_n}$
if and only if 
so does $P\rho P^{-1}$. 
We also say that 
a generalized ${\mathcal B}_n$-representation with 
Borel mold  
$\{ (U_i, \rho_i ) \}_{i \in I}$ 
{\it satisfies the $(\ast)$-condition} 
with respect to 
$\{ \alpha_{i, j} \}
_{ (i, j) \in {\mathcal I}_n}$ 
if the same condition holds.  

\bigskip

For an ${\mathcal I}_n$-indexed set $A = \{ \alpha_{i, j} \}
_{ (i, j) \in {\mathcal I}_n}$,  
the subfunctor ${\mathcal EqBB}_n(\Gamma)_{A}$ 
of ${\mathcal EqBB}_n(\Gamma)$ is 
defined by 
${\mathcal EqBB}_n(\Gamma)_A(X) := 
\{ 
 \rho \in 
 {\mathcal EqBB}_n(\Gamma)(X) \mid 
 \rho \mbox{ satisfies the $(\ast)$-condition with respect to } 
 A 
 \}$
for a scheme $X$.  
Let $k$ be a field and let $\rho \in {\mathcal EqBB}_n(\Gamma)(k)$. 
Let $\tilde{\rho} : \Gamma 
\to {\mathcal B}_n(k)$ 
be a representative of $\rho$.   
Then it is  easy to check that 
$\rho$ satisfies the $(\ast)$-condition with
respect to some ${\mathcal I}_n$-indexed subset 
$A = \{ \alpha_{i, j} \}
_{ (i, j) \in {\mathcal I}_n}$ of $\Gamma$. 
Hence we have 
${\mathcal EqBB}_n(\Gamma)(k) = \cup_{A} \; 
{\mathcal EqBB}_n(\Gamma)_A(k)$, 
where the union runs the ${\mathcal I}_n$-indexed subsets of $\Gamma$.   

In the sequel, we show that ${\mathcal EqBB}_n(\Gamma)_A$
is an affine scheme for each ${\mathcal I}_n$-indexed subset $A$
of $\Gamma$. 
Let us define the subfunctor ${\rm B}_n(\Gamma)_{B, A}$ of 
${\rm B}_n(\Gamma)$ by 
\[
{\rm B}_n(\Gamma)_{B, A}(X) := \{ 
\sigma \in {\rm B}_n(\Gamma)(X) \mid \sigma \mbox{ satisfies the 
$(\ast)$-condition for }A \}
\]
for a scheme $X$.  
We easily verify that 
${\rm B}_n(\Gamma)_{B, A}$ is an affine open subscheme of 
${\rm B}_n(\Gamma)$. 
The action of ${\rm B}_n$ on ${\rm B}_n(\Gamma)$ by 
$\rho \mapsto Q \rho Q^{-1}$ induces the one of ${\rm B}_n$ 
on ${\rm B}_n(\Gamma)_{B, A}$. 
There exists a canonical morphism ${\rm B}_n(\Gamma)_{B, A} \to 
{\mathcal EqBB}_n(\Gamma)_A$. Then we obtain the following:

\begin{lemma}\label{lemma-eqbb}
	The morphism ${\rm B}_n(\Gamma)_{B, A} \to 
{\mathcal EqBB}_n(\Gamma)_A$ is a ${\rm B}_n$-principal 
fiber bundle. 
In particular, the functor ${\mathcal EqBB}_n(\Gamma)_A$
is an affine scheme. 
\end{lemma}

As a corollary of Lemma \ref{lemma-eqbb}, we have: 

\begin{corollary}\label{cor-eqbb}
	The functor ${\mathcal EqBB}_n(\Gamma)$ is representable. 
\end{corollary}

\prf The statement follows from that 
the sheaf ${\mathcal EqBB}_n(\Gamma)$ is 
covered by affine schemes ${\mathcal EqBB}_n(\Gamma)_{A}$.
\qed

\bigskip

Before proving Lemma \ref{lemma-eqbb}, we need several 
preparations and long discussions. 
Let $\rho$ be the universal ${\mathcal B}_n$-representation 
on ${\rm B}_n(\Gamma)$. 
Fix an ${\mathcal I}_n$-indexed subset 
$A = \{ \alpha_{i, j} \}_{ (i, j) \in {\mathcal I}_n}$ of $\Gamma$. 
We denote the coordinate ring of the affine scheme 
${\rm B}_n(\Gamma)_{B, A}$ by $R$. 
We define $\eta_{ij}(\gamma) \in J_{ij}(R)$ 
and $\varepsilon_{ij}(\gamma) \in R$ 
for $\gamma \in \Gamma$ and $(i, j) \in {\mathcal I}_n$ by induction.
First we define $\varepsilon_{11}(\gamma) := \rho(\gamma)_{11}/
\rho(\alpha_{11})_{11}$ and $\eta_{11}(\gamma) := 
\rho(\gamma) - \varepsilon_{11}(\gamma)\rho(\alpha_{11})$. 
Suppose that $\eta_{i'j'}(\gamma) \in J_{i'j'}(R)$ 
and $\varepsilon_{i'j'}(\gamma) \in R$ are defined for 
each $\gamma \in \Gamma$ and for each $(i', j') < (i, j)$. 
Now we define $\varepsilon_{ij}(\gamma) := (\eta_{i'j'}(\gamma))_{ij}/
(\eta_{i'j'}(\alpha_{ij}))_{ij}$ and $\eta_{ij}(\gamma) := 
\eta_{i'j'}(\gamma) - \varepsilon_{ij}(\gamma)\eta_{i'j'}(\alpha_{ij})$, 
where $(i', j')$ is the previous index of $(i, j)$. 
Here we remark that $(\eta_{i', j'}(\alpha_{ij}))_{ij}
\in R^{\times}$. 
Set $\tau_{11}:=(\rho_{11}(\alpha_{11}))_{11} \in R^{\times}$ and  
$\tau_{ij}:=(\eta_{i', j'}(\alpha_{ij}))_{ij}
\in R^{\times}$ for $(i, j) \in {\mathcal I}_n \setminus \{ (1, 1) \}$. 

For $Q = (b_{ij}) \in {\rm B}_n$, let us denote by $Q^{\ast} : R \to R$ the 
isomorphism associated to the adjoint action on ${\rm B}_n(\Gamma)_{B, A}$ 
by $\rho \mapsto Q\rho Q^{-1}$.  We also denote by 
$Q^{\ast} : {\mathcal B}_n(R) \to {\mathcal B}_n(R)$ the isomorphism 
defined by $X = (x_{ij}) \mapsto Q^{\ast}(X) := (Q^{\ast}(x_{ij}))$.   
Under the action of ${\rm B}_n$ on $R$, 
the elements $\tau_{ij}$ are semi-invariant. 
Indeed, for $Q = (b_{ij}) \in {\rm B}_n$ we have 
$Q^{\ast}\rho(\gamma)=Q \rho(\gamma)Q^{-1}$, 
$Q^{\ast}\eta_{ij}(\gamma)=Q\eta_{ij}(\gamma)Q^{-1}$, 
and $Q^{\ast} \tau_{ij} = ({b_{ii}}/{b_{jj}}) \tau_{ij}$. 
Hence $\tau_{ij}\tau_{jk}\tau_{ik}^{-1}$ is 
${\rm B}_n$-invariant. 
We also see that $\varepsilon_{ij}(\gamma)$ is 
${\rm B}_n$-invariant. 

Let us introduce the ${\rm B}_n$-invariant subalgebra of 
$R$.   

%
%

\begin{definition}\rm
Let $R^{\rm ch}$ be the ${\rm B}_n$-invariant subalgebra of $R$. 
More precisely, we define $R^{\rm ch} := 
\{ x \in R \mid \sigma^{\ast}(x) = 1 \otimes x \}$, where 
$\sigma^{\ast} : R \to A({\rm B}_n)\otimes_{{\Bbb Z}} R$ is the ring homomorphism 
associated to the group action $\sigma : {\rm B}_n \times {\rm B}_n(\Gamma)_{B, A} 
\to {\rm B}_n(\Gamma)_{B, A}$.  
Here we denote by $A({\rm B}_n)$ the coordinate ring of ${\rm B}_n$.    
We define the affine scheme 
${\rm Ch}_n(\Gamma)_{B, A} := {\rm Spec}R^{\rm ch}$.   
\end{definition}

The next lemma is a key for proving Lemma \ref{lemma-eqbb}: 

\begin{lemma}\label{lemma:conj}
There exists an upper triangular invertible matrix $Q \in \tilde{B}_n(R) := 
\{ (b_{ij}) \in {\rm GL}_n(R) \mid b_{ij}=0 \mbox{ for } i>j \}$ 
such that all entries of $Q\rho(\gamma) Q^{-1}$ are contained in $R^{\rm ch}$ 
for each $\gamma \in \Gamma$. 
\end{lemma}

From now on, we concentrate ourselves into proving Lemma \ref{lemma:conj}. 
Note that the representation $\rho$ satisfies 
\[
\rho(\gamma)=\sum_{(i, j) \in {\mathcal I}_n} 
\varepsilon_{i, j}(\gamma)\eta_{i'j'}(\alpha_{ij}) \] 
for each $\gamma \in \Gamma$. 
Here we put $\eta_{1'1'}(\alpha_{11}) = \rho(\alpha_{11})$. 
For proving Lemma \ref{lemma:conj}, we only need to show that 
there exists $Q \in \tilde{B}_n(R)$ such that 
$Q\eta_{i'j'}(\alpha_{ij})Q^{-1} \in {\mathcal B}_n(R^{\rm ch})$ 
for each $(i, j) \in {\mathcal I}_n$.

\begin{definition}\rm
For $X \in {\mathcal B}_n(R)$, we say that $X$ is {\it canonical} 
if the action of ${\rm B}_n$ on (entries of) $X$ is described by 
$Q^{\ast}X = QXQ^{-1}$ for $Q \in {\rm B}_n$. 
Note that $\rho(\gamma)$ and $\eta_{ij}(\gamma)$ are canonical for each 
$\gamma \in \Gamma$ and $(i, j) \in {\mathcal I}_n$. 
\end{definition}

\begin{definition}\rm
For $(i, j) \in {\mathcal I}_n$, we define the {\it convex hull } 
$\overline{(i, j)}$ of $(i, j)$ as the subset 
$\overline{(i, j)} := \{ (k , \ell) \in {\mathcal I}_n 
\mid k \le i,\; \ell \ge j \}$ of ${\mathcal I}_n$.  
We also define the {\it convex hull } 
$\overline{S}$ of a subset $S$ of ${\mathcal I}_n$ 
by $\displaystyle \overline{S} := \cup_{(i, j) \in S} \: \overline{(i, j)}$.  
\end{definition}

\begin{definition}\rm 
Let $X = (x_{ij})\in {\mathcal B}_n(R)$. 
We define the {\it support} of $X$ by 
${\rm Supp}X := \{ (i, j) \in {\mathcal I}_n \mid x_{ij} \neq 0 \}$. 
We say that $X$ is {\it $(i, j)$-shaped} if 
$x_{ij} \in R^{\times}$ and ${\rm Supp}X \subseteq \overline{(i, j)}$. 
We also say that $X$ is {\it well-shaped} if 
$X = 0$ or $X$ is $(i, j)$-shaped for some $(i, j) \in {\mathcal I}_n$. 
\end{definition}

Set $Y(i,j) := \eta_{i'j'}(\alpha_{ij})$ for $(i, j) \in {\mathcal I}_n$. 
(Recall $Y(1, 1) := \eta_{1'1'}(\alpha_{11}) = \rho(\alpha_{11})$.) 
Note that $Y(i,j)$ is canonical and that the $(i, j)$-entry of $Y(i,j)$ is 
contained in $R^{\times}$. 
Now we show the following lemma. 

\begin{lemma}
For each $(i, j) \in {\mathcal I}_n$, there exists an $(i, j)$-shaped canonical matrix 
$X(i, j) \in {\mathcal B}_n(R)$ such that 
\[ \displaystyle Y(i, j) = X(i, j) + 
\sum_{(k, \ell) \in \overline{{\rm Supp}Y(i, j)} \setminus \overline{(i, j)}} a_{k, \ell}(i, j) X(k, \ell), 
\]  
where $a_{k, \ell}(i, j)$ is ${\rm B}_n$-invariant. 
\end{lemma}

\prf
Since $Y(1, n)$ is a $(1, n)$-shaped canonical matrix, we set $X(1, n):=Y(1, n)$. 
Suppose that we can define a $(k, \ell)$-shaped canonical matrix 
$X(k, \ell)$ for $(k, \ell) > (i, j)$ and that 
the equalities above hold. 
Let us consider the case $(i, j)$. 
Set $J := \overline{{\rm Supp}Y(i, j)} \setminus \overline{(i, j)}$. 
If $J = \emptyset$, then set $X(i, j) := Y(i, j)$. 
Assume that $J \neq \emptyset$. 
Let $(s, t)$ be the minimum element of $J$. 
Then $(s, t) \in {\rm Supp}Y(i, j)$. 
If $(s', t') \neq (s, t)$ for $(s', t') \in {\rm Supp}Y(i, j)$, then 
$(s, t) \not\in \overline{(s', t')}$.   
Then the $(s, t)$-entry of $B Y(i, j) B^{-1}$ is equal to $y b_{ss}/b_{tt}$, 
where $B=(b_{**})$ and $y$ is the $(s, t)$-entry of $Y(i, j)$. 
In other words, $y$ is semi-invariant. 
Remark that $(s, t)>(i, j)$. 
The $(s, t)$-entry $y'$ of $X(s, t)$ is a unit, and 
$a_{st}(i, j) := y/y'$ is ${\rm B}_n$-invariant. 
The new matrix $Y':=Y(i, j)-a_{st}(i, j)X(s, t)$ is a canonical matrix and 
${\rm Supp}Y' \subseteq \overline{(i, j)} \cup (J \setminus \{ (s, t) \})$. 
Instead of $Y(i, j)$ and $J$, we consider $Y'$ and $J' := \overline{{\rm Supp}Y'} 
\setminus \overline{(i, j)}$. 
Then $J' = \emptyset$ or the minimal element of $J'$ is bigger than $(s, t)$. 
By induction on the minimum elements of $J'$, we can obtain 
${\rm B}_n$-invariants $a_{k, \ell}(i, j)$ for $(k, \ell) \in J$ such that 
$X(i, j) := Y(i, j) - \sum_{(k, \ell)} a_{k\ell}(i, j) X(k, \ell)$ is an 
$(i, j)$-shaped canonical matrix.   
By repeating this discussion, we can obtain an $(i, j)$-shaped canonical 
matrix $X(i, j)$ for each $(i, j) \in {\mathcal I}_n$. 
\qed

\bigskip 

From the lemma above, we obtained well-shaped canonical matrices $X(i, j)$ 
from $Y(i, j) = \eta_{i'j'}(\alpha_{ij})$. 
For proving Lemma \ref{lemma:conj}, we only need to verify that 
there exists $Q \in \tilde{B}_n(R)$ such that $QX(i, j)Q^{-1} \in 
{\mathcal B}_n(R^{\rm ch})$. 
The next lemma is useful for the discussion below. 

\begin{lemma}\label{lemma:red} 
Let $X$ be a canonical matrix of ${\mathcal B}_n(R)$. Set $J := \overline{{\rm Supp}X}$. 
Then there exist ${\rm B}_n$-invariants $a_{ij}$ for $(i, j) \in J$ such that 
\[
X = \sum_{(i, j) \in J} a_{ij} X(i, j). 
\]
\end{lemma}

\prf
If $J=\emptyset$, then the statement is trivial. 
Suppose that $J \neq \emptyset$. 
Let $(i, j) \in J$ be the minimum element of $J$. 
The $(i, j)$-entry $x_{ij}$ of $X$ is semi-invariant, and hence 
$a_{ij} := x_{ij}/X(i, j)_{ij}$ is ${\rm B}_n$-invariant. Here
we denote by $X(i,j)_{ij}$ the $(i, j)$-entry of $X(i, j)$. 
The new matrix $X' := X - a_{ij}X_{ij}$ is a canonical matrix. 
If $X' \neq 0$, then $J' := \overline{{\rm Supp}X'} \subset J$ and  
the minimum element of $J'$  
is bigger than $(i, j)$. 
By induction on the minimum of $J = \overline{{\rm Supp}X}$, we can prove the
statement. 
\qed

\begin{definition}\label{def:r'}\rm
Let $R'$ be the subalgebra of $R$ over $R^{\rm ch}$ generated by 
the following 
elements:  
\[
C := \left\{ 
\begin{array}{cl}  
(X(1, i)_{1i})^{\pm 1} & ([\frac{n+1}{2}]+1 \le i \le n), \\ 
(X(i, n)_{in})^{\pm 1} & (2 \le i \le [\frac{n+1}{2}]), \\
X(1, 1)_{1i} & (2 \le i \le n), \\
X(i, i)_{ji} & (3 \le i \le n, \; 2 \le j \le i-1)  \\
\end{array} 
\right\}.  
\]
\end{definition}

\begin{lemma}\label{lemma:xij}
For each $(i, j) \in {\mathcal I}_n$, $X(i, j) \in {\mathcal B}_n(R')$. 
\end{lemma} 

\prf
Note that $X(i, i)_{ii} \in (R^{\rm ch})^{\times}$ for $1 \le i \le n$.    
It is easy to see that $X(1, 1), X(1, n) \in {\mathcal B}_n(R')$. 
First, we verify that $X(n, n) \in {\mathcal B}_n(R')$. 
Remark that $X(n,n)_{in} \in C$ for $2 \le i \le n-1$ and 
$X(n, n)_{nn} \in (R^{\rm ch})^{\times}$. 
For proving that $X(n, n)_{1n} \in R'$, 
let us consider the canonical matrix $X(1, 1)X(n, n)$. 
By Lemma \ref{lemma:red}, we see that 
$X(1, 1)X(n, n) = a(1, n)X(1, n)$ for some $a(1, n) \in R^{\rm ch}$ 
because ${\rm Supp} (X(1, 1)X(n, n)) \subseteq \{(1, n) \}$. 
Comparing the $(1, n)$-entries, we have 
\[
X(1, 1)_{11}X(n, n)_{1n}+ \sum_{k=2}^n X(1, 1)_{1k}X(n, n)_{kn} = a(1, n)X(1, n)_{1n}.
\]
Hence we see that 
\[
X(n, n)_{1n} = X(1, 1)_{11}^{-1}( a(1, n)X(1, n)_{1n} - \sum_{k=2}^n X(1, 1)_{1k}X(n, n)_{kn})  
\] 
and that $X(n, n)_{1n} \in R'$. 
Therefore we have $X(n, n) \in {\mathcal B}_n(R')$.  

Next, we show that $X(1, i)_{1i} \in (R')^{\times}$ for $2 \le i \le n$ and 
that $X(i, n)_{in} \in (R')^{\times}$ for $2 \le i \le n$.  
Let us consider the canonical matrix $X(1, i)X(i, n)$. 
The support is contained in $\{ (1, n) \}$, and 
$X(1, i)X(i, n) \\ = a X(1, n)$ for some $a \in R^{\rm ch}$. 
Comparing the $(1, n)$-entries, we have 
\[
X(1, i)_{1i}X(i, n)_{in} = a X(1, n)_{1n}.
\] 
Since $X(1, i)_{1i}, X(i, n)_{in}, X(1, n)_{1n} \in R^{\times}$, 
we have $a \in (R^{\rm ch})^{\times}$. 
Because $(X(1, i)_{1i})^{\pm 1} \in C$ for $[\frac{n+1}{2}]+1 \le i \le n$ and   
$(X(i, n)_{in})^{\pm 1} \in C$ for $2 \le i \le [\frac{n+1}{2}]$, 
we see that $X(1, i)_{1i} \in (R')^{\times}$ for $2 \le i \le n$ and 
that $X(i, n)_{in} \in (R')^{\times}$ for $2 \le i \le n$.

Third, we show that $X(1, i) \in {\mathcal B}_n(R')$ for $1 < i < n$. 
We have known that $X(1, n) \in {\mathcal B}_n(R')$. 
Let us consider the canonical matrices $X(1, i)X(j, j)$ for $j > i$. 
By Lemma \ref{lemma:red} there exist ${\rm B}_n$-invariants $a(1, k)$  
for $j \le k \le n$ such that 
$X(1, i)X(j, j) = a(1, j)X(1, j)+ \sum_{k > j} a(1, k)X(1, k)$ 
since ${\rm Supp}X(1, i)X(j, j) \subseteq \overline{\{ (1, j) \}}$.  
Comparing the $(1, j)$-entries, we have 
\[
\begin{array}{r}
X(1, i)_{1i}X(j, j)_{ij}+X(1, i)_{1, i+1}X(j, j)_{i+1, j}+ 
\cdots + X(1, i)_{1j}X(j, j)_{jj} \\ 
= a(1, j)X(1, j)_{1j}. 
\end{array}
\]
We have seen that $X(1, i)_{1i} \in R'$. 
Assume that $X(1, i)_{1k} \in R'$ 
for $i \le k  \le j-1$.  
By the equality above, we have 
\[
\begin{array}{lr}
X(1, i)_{1j} = & X(j, j)_{jj}^{-1}(a(1, j)X(1, j)_{1j} - X(1, i)_{1i}X(j, j)_{ij} - \cdots  \\ 
 & - X(1, i)_{1, j-1}X(j, j)_{j-1, j} ). 
\end{array}
\]
Since $X(j, j)_{ij}, \ldots, X(j, j)_{j-1, j} \in C$ and 
$X(1, j)_{1j} \in R'$, we obtain $X(1, i)_{1j} \in R'$. 
By induction on $j$, we have $X(1, i)_{1j} \in R'$ for $i \le j \le n$. 
Thus $X(1, i) \in {\mathcal B}_n(R')$. 


Finally, we show that $X(i, j) \in {\mathcal B}_n(R')$ for each $(i, j) \in {\mathcal I}_n$. 
If $i=1$, then we have checked it. 
Hence we may assume that $i>1$. 
By Lemma \ref{lemma:red} we see that 
\[
X(1, i)X(i, j) = a(1, j)X(1, j)+ \cdots +a(1, n)X(1, n) 
\]
for suitable  ${\rm B}_n$-invariants $a(1, k) \; (j \le k \le n)$. 
In particular, we obtain $X(1, i)X(i, j) \in {\mathcal B}_n(R')$.    
The $(1, k)$-entry of $X(1, i)X(i, j)$ is $X(1, i)_{1i}X(i, j)_{ik}$ 
for $j \le k \le n$. 
Since $X(1, i)_{1i} \in (R')^{\times}$, $X(i, j)_{ik} \in R'$.
Suppose that there exists $1 \le \ell \le i-1$ such that $X(i, j)_{mk} \in R'$ for $m>\ell$ and $j \le k \le n$. 
Then we prove that $X(i, j)_{\ell k} \in R'$ for $j \le k \le n$.  
By using Lemma \ref{lemma:red} again, we see that 
\[
X(1, \ell)X(i, j) = a(1, j)X(1, j)+ \cdots +a(1, n)X(1, n) \in {\mathcal B}_n(R') 
\]
for suitable  ${\rm B}_n$-invariants $a(1, k) \; (j \le k \le n)$. 
For $j \le k \le n$, the $(1, k)$-entry of $X(1, \ell)X(i, j)$ is 
\[
X(1, \ell)_{1\ell}X(i, j)_{\ell k} + X(1, \ell)_{1,\ell+1}X(i, j)_{\ell+1, k} + \cdots 
+ X(1, \ell)_{1i}X(i, j)_{i k} \in R'.
\]
Since $X(1, \ell)_{1\ell} \in (R')^{\times}$ and $X(1, \ell) \in {\mathcal B}_n(R')$, 
we have $X(i, j)_{\ell k} \in R'$ for $j \le k \le n$ by the hypothesis.   
By induction on $\ell$, we see that $X(i, j)_{\ell k} \in R'$ for $1 \le \ell \le i$ and $j \le k \le n$.
Therefore we proved that $X(i, j) \in {\mathcal B}_n(R')$ for each $(i, j) \in {\mathcal I}_n$. 
\qed

\bigskip 

Now we can prove Lemma \ref{lemma:conj}. 

\bigskip 

{\it Proof of Lemma \ref{lemma:conj}}. 
By the long discussion above, we only need to show that 
there exists $Q \in \tilde{B}_n(R)$ such that $QX(i, j)Q^{-1} \in 
{\mathcal B}_n(R^{\rm ch})$. 
For $\{ X(i, j) \mid (i, j) \in {\mathcal I}_n \}$ we introduced 
the set $C$ in Definition \ref{def:r'}. 
Let $Q \in \tilde{B}_n(R)$. 
For $\{ QX(i, j)Q^{-1} \mid (i, j) \in {\mathcal I}_n \}$ we define 
the set $Q^{\ast}C$ in a similar way. 
From now on, we prove that there exists $Q \in \tilde{B}_n(R)$ 
such that any element of $Q^{\ast}C$ is $0$ or $1$. 
More precisely, we prove that 
\begin{eqnarray*}
(QX(1, i)Q^{-1})_{1i} & = & 1 \hspace{2ex} ( [\frac{n+1}{2}]+1 \le i \le n) \\ 
(QX(i, n)Q^{-1})_{in} & = & 1 \hspace{2ex} ( 2 \le i \le [\frac{n+1}{2}]) \\
(QX(1, 1)Q^{-1})_{1i} & = & 0 \hspace{2ex} ( 2 \le i \le n) \\  
(QX(i, i)Q^{-1})_{ji} & = & 0 \hspace{2ex} ( 3 \le i \le n, 2 \le j \le i-1). \\ 
\end{eqnarray*}

Let us find $Q = (q_{ij}) \in \tilde{B}_n(R)$. 
First, we set  
\[
(q_{11}, q_{22}, \ldots, q_{nn}) 
\displaystyle := (1, \frac{\tau_{1, n}}{\tau_{2, n}}, 
\frac{\tau_{1, n}}{\tau_{3, n}}, 
\ldots, \frac{\tau_{1, n}}{\tau_{[\frac{n+1}{2}], n}}, \tau_{1, [\frac{n+1}{2}]+1}, 
\ldots, \tau_{1, n-1}, \tau_{1, n}).
\] 
Here recall that $\tau_{ij} = \eta_{i'j'}(\alpha_{ij})_{ij} = Y(i, j)_{ij} 
=X(i, j)_{ij}$. 
Then it is easy to see that  
$(QX(1, i)Q^{-1})_{1i} = 1$ for $[\frac{n+1}{2}]+1 \le i \le n$ and  
$(QX(i, n)Q^{-1})_{in} = 1$ for $2 \le i \le [\frac{n+1}{2}]$.  

Next, let us define $q_{ij} \in R$ for $i<j$. 
Set $Q^{-1} = (q'_{\ast\ast})$. 
Let $X = (x_{\ast\ast}) \in {\mathcal B}_n(R)$. 
From $QQ^{-1} = I_n$, we have $q'_{kk} = q_{kk}^{-1}$ and  
\[
q_{ij}q'_{jj}+q_{i, j-1}q'_{j-1, j} + \cdots + q_{i, i+1}q'_{i+1, j} + q_{ii}q'_{ij} = 0 
\]
for $i<j$. 
Since 
\[ 
q'_{ij} = -q_{ii}^{-1}(q_{ij}q'_{jj} + q_{i, j-1}q'_{j-1, j} + \cdots + q_{i, i+1}q'_{i+1, j}), 
\]
we see that $q'_{ij}$ can be expressed by 
$\{ q_{k\ell} \mid (k, \ell) \le (i, j) \}$.  
The $(i, j)$-entry of $QXQ^{-1}$ is 
\begin{eqnarray*}
& & q_{ij}x_{jj}q'_{jj} + q_{ii}x_{ii}q'_{ij} + 
\sum_{(k, \ell) \neq (i, i), (j, j)} q_{ik}x_{k\ell}q'_{\ell j} \\
& = & (x_{jj} - x_{ii})q_{jj}^{-1} q_{ij} + (\text{ rational function of } 
\{ q_{k\ell} \mid (k, \ell) < (i, j) \} ). 
\end{eqnarray*} 
Assume that $q_{k\ell} \in R$ is defined for each $(k, \ell) < (i, j)$.  
Now we define $q_{ij} \in R$.
If $i=1$, then put $X = X(1, 1)$. 
Since $x_{11} \in R^{\times}$ and $x_{jj} = 0$, 
the $(1, j)$-entry of $QX(1, 1)Q^{-1}$ is 
$-x_{11}q_{jj}^{-1}q_{ij} + ( \text{lower term} )$, and hence 
we can find $q_{ij}$ satisfying the equation 
$(QX(1, 1)Q^{-1})_{1j}=0$.   
If $i>1$, then put $X = X(j, j)$. 
Since $x_{ii}=0$ and $x_{jj} \in R^{\times}$, 
the $(i, j)$-entry of $QX(j, j)Q^{-1}$ is 
$x_{jj}q_{jj}^{-1}q_{ij} + ( \text{lower term} )$, and hence 
we can find $q_{ij}$ satisfying the equation 
$(QX(j, j)Q^{-1})_{ij}=0$. 
By induction on $(i, j) \in {\mathcal I}_n$, we can 
define $Q = (q_{ij})$ satisfying the equations.   
In particular, any element of $Q^{\ast}C$ is $0$ or $1$. 

By Lemma \ref{lemma:xij} we have $X(i, j) \in {\mathcal B}_n(R')$. 
Similarly, we see that all entries of $QX(i, j)Q^{-1}$ are 
contained in the algebra generated by $Q^{\ast}C$ over $R^{\rm ch}$. 
Since any element of $Q^{\ast}C$ is $0$ or $1$, Q$X(i, j)Q^{-1} \in 
{\mathcal B}_n(R^{\rm ch})$ for $(i, j) \in {\mathcal I}_n$. 
\qed 

\bigskip

\medskip

Now we can finish the proof of Lemma \ref{lemma-eqbb}.

\bigskip 

{\it Proof of Lemma \ref{lemma-eqbb}.}  
By Lemma \ref{lemma:conj} there exists $Q \in \tilde{B}_n(R)$ 
such that $Q\rho(\gamma)Q^{-1} \in {\mathcal B}_n(R^{\rm ch})$ 
for each $\gamma \in \Gamma$. 
The inclusion $R^{\rm ch} \to R$ induces the morphism 
$\pi : {\rm B}_n(\Gamma)_{B, A} \to {\rm Ch}_n(\Gamma)_{B, A}$. 
We have the section $s : {\rm Ch}_n(\Gamma)_{B, A} 
\to {\rm B}_n(\Gamma)_{B, A}$ associated to  
the ${\mathcal B}_n$-representation $\rho':= Q\rho Q^{-1} : \Gamma \to 
{\mathcal B}_n(R^{\rm ch})$ with Borel mold.

We show that the morphism $\phi : {\rm B}_n 
\times {\rm Ch}_n(\Gamma)_{B, A}  \to 
{\rm B}_n(\Gamma)_{B, A}$ 
associated to the ${\mathcal B}_n$-representation 
$\tilde{Q}\rho' \tilde{Q}^{-1}$ gives 
an isomorphism.
Here we denote by $\tilde{Q}$ the universal matrix of ${\rm B}_n$. 
For a scheme $Z$, the morphism $\phi_{\ast}(Z) : h_{{\rm B}_n}(Z) 
\times h_{{\rm Ch}_n(\Gamma)_{B, A}}(Z)  
 \to h_{{\rm B}_n(\Gamma)_{B, A}}(Z)$ is injective because of 
Proposition \ref{prop:actionfree}. 
Let us prove that $\phi_{\ast}(Z)$ is surjective. 
We denote by $\overline{Q} \in {\rm B}_n(R)$ the image of 
$Q$ by $\tilde{B}_n(R) \to {\rm B}_n(R)$.  
The representation $\rho' = \overline{Q}\rho \overline{Q}^{-1}$ corresponds to 
${\rm B}_n(\Gamma)_{B, A} \stackrel{\pi}{\to} {\rm Ch}_n(\Gamma)_{B, A} 
\stackrel{s}{\to}  {\rm B}_n(\Gamma)_{B, A}$. 
Then ${\rm B}_n(\Gamma)_{B, A} \stackrel{(\overline{Q}^{-1}, \rho')}{\to} 
{\rm B}_n \times {\rm B}_n(\Gamma)_{B, A} \stackrel{\sigma}{\to} 
{\rm B}_n(\Gamma)_{B, A}$ is the identity. 
The morphism ${\rm B}_n(\Gamma)_{B, A} \stackrel{(\overline{Q}^{-1}, \rho')}{\to} 
{\rm B}_n \times {\rm B}_n(\Gamma)_{B, A}$ factors through 
${\rm B}_n(\Gamma)_{B, A} \stackrel{(\overline{Q}^{-1}, \pi)}{\to} 
{\rm B}_n \times {\rm Ch}_n(\Gamma)_{B, A} 
\stackrel{(id, s)}{\to} {\rm B}_n \times {\rm B}_n(\Gamma)_{B, A}$. 
Hence $\phi_{\ast}(Z)$ is surjective. 
Therefore $\phi$ is an isomorphism.

The morphism 
${\rm B}_n(\Gamma)_{B, A} \to {\rm Ch}_n(\Gamma)_{B, A}$
gives a ${\rm B}_n$-principal fiber bundle.
We can check that the functor ${\mathcal EqBB}_n(\Gamma)_{A}$ is 
representable by 
${\rm Ch}_n(\Gamma)_{B, A}$. 
This completes the proof of Lemma 
\ref{lemma-eqbb}. 
\qed

\bigskip

By Corollary \ref{cor-eqbb}, we see that 
${\mathcal EqBB}_n(\Gamma)$ is representable. 
We introduce the following definition.

\begin{definition}\rm
The scheme ${\rm Ch}_n(\Gamma)_B$ which represents 
the functor ${\mathcal EqBB}_n(\Gamma) = {\mathcal EqB}_n(\Gamma)$
is called the {\it moduli of representations 
with Borel mold} of degree $n$ for $\Gamma$. 
It is also called the {\it character variety with Borel mold} 
of degree $n$ for $\Gamma$. 
\end{definition}

\begin{remark}\rm
The canonical morphism 
$\pi : {\rm Rep}_n(\Gamma)_B \to {\rm Ch}_n(\Gamma)_B$ 
is a principal fiber bundle with fiber ${\rm PGL}_n$. 
The canonical morphism $\pi' : 
{\rm B}_n(\Gamma)_B \to {\rm Ch}_n(\Gamma)_B$ 
is also a principal fiber bundle 
with fiber ${\rm B}_n$. 
These principal fiber bundles have a local trivialization 
with respect to Zariski topology.  
They are universal geometric quotients in \cite{GIT}.  

The construction of the moduli of representations with Borel 
mold gives us the following diagram:
\[
\begin{array}{ccc}
{\rm B}_n(\Gamma)_{B}\times {\rm PGL}_n & \stackrel{f}{\to} 
& {\rm Rep}_n(\Gamma)_B \\
\downarrow {p_1} & & \downarrow {\pi} \\
{\rm B}_n(\Gamma)_B & \stackrel{\pi'}{\to} & {\rm Ch}_n(\Gamma)_B, 
\end{array}
\]
where 
\[
\begin{array}{ccccc}
f & : & {\rm B}_n(\Gamma)_{B}\times {\rm PGL}_n & \to 
& {\rm Rep}_n(\Gamma)_B \\
& & (\rho, Q) & \mapsto & Q^{-1}\rho Q  
\end{array}
\]
and $p_1 : {\rm B}_n(\Gamma)_{B}\times {\rm PGL}_n 
\to {\rm B}_n(\Gamma)_{B}$ is the first projection. 
The morphism $f$ is a principal fiber bundle with fiber 
${\rm B}_n$ which has a local trivialization 
with respect to Zariski topology.  
\end{remark}

\bigskip


\begin{lemma}\label{lemma:finitetype}
Let $\Gamma$ be a finitely generated group or monoid.  
Then ${\rm Ch}_n(\Gamma)_B$ is of finite type over 
${\Bbb Z}$. 
\end{lemma}

\prf
Since the representation variety ${\rm Rep}_n(\Gamma)$ 
is of finite type over ${\Bbb Z}$ when $\Gamma$ 
is finitely generated, so is 
a subscheme 
${\rm Rep}_n(\Gamma)_B$.  
The principal fiber bundle 
$\pi : {\rm Rep}_n(\Gamma)_B \to {\rm Ch}_n(\Gamma)_B$
with fiber ${\rm PGL}_n$ has a local 
trivialization with respect to 
Zariski topology, and hence 
${\rm Ch}_n(\Gamma)_B$ is of finite type over ${\Bbb Z}$.  
\qed

\bigskip 

\begin{remark}\rm
Lemma \ref{lemma:finitetype} can be verified by 
investigating the invariants directly. 
Since ${\rm Rep}_n(\Gamma)_B$ is quasi-compact, ${\rm Ch}_n(\Gamma)_B$ is also 
quasi-compact. It is essential to prove that 
${\rm Ch}_n(\Gamma)_B$ is locally of finite type over ${\mathbb Z}$. 
Let $R^{\rm ch}$ be the affine ring of ${\rm Ch}_n(\Gamma)_{B, A}$.  
Let $\rho' : \Gamma \to {\mathcal B}_n(R^{\rm ch})$ be the 
${\mathcal B}_n$-representation with Borel mold in the proof of 
Lemma \ref{lemma-eqbb}. 
For a generator $\{ \alpha_{i} \}_{i=1}^{N}$ of $\Gamma$, 
we consider the set $S$ of all entries of $\rho'(\alpha_i)$ for $i=1, 2, \ldots, N$.   
Then $R^{\rm ch}$ is generated by $S$ over ${\mathbb Z}$. 
Indeed, let $R_{0}$ be the subalgebra of $R^{\rm ch}$ generated by $S$ over ${\mathbb Z}$. 
We can define the ${\mathcal B}_n$-representation $\rho'' : \Gamma 
\to {\mathcal B}_n(R_0)$ with Borel mold such that 
$\rho'' \otimes_{R_0} R^{\rm ch} = \rho'$. 
We define the section of ${\rm B}_n(\Gamma)_{B, A} \to 
{\rm Ch}_n(\Gamma)_{B, A} \to {\rm Spec}R_0$ associated to $\rho''$. 
We can show that for $f \in {\rm Hom}({\rm Ch}_n(\Gamma)_{B, A}, {\mathbb A}^1_{\mathbb Z})$ 
there exists a unique $f' \in {\rm Hom}({\rm Spec}R_0, {\mathbb A}^1_{\mathbb Z})$ 
such that ${\rm Ch}_n(\Gamma)_{B, A} \to {\rm Spec}R_0 \stackrel{f'}{\to} 
{\mathbb A}^1_{\mathbb Z}$ is $f$. Hence we have $R^{\rm ch} = R_0$.   
Therefore $R^{\rm ch}$ is finitely generated over ${\mathbb Z}$. 
\end{remark} 

\bigskip 

\begin{proposition}\label{prop-sep}
The moduli ${\rm Ch}_n(\Gamma)_B$ is 
a separated scheme over ${\Bbb Z}$.
\end{proposition}

\bigskip

Before proving Proposition \ref{prop-sep}, we introduce the following 
lemma. 
We define the subgroup scheme $\tilde{B}_n$ of ${\rm GL}_n$ 
by $\tilde{B}_n := \{ (b_{ij}) \in {\rm GL}_n \mid b_{ij}=0 \mbox{ for each }
i>j \}$. 

\begin{lemma}\label{lemma-sep}
Let $R$ be a valuation ring and $K$ its quotient field. 
Suppose that $P \in \tilde{B}_n(K)$ satisfies $P {\mathcal B}_n(R) P^{-1} 
={\mathcal B}_n(R)$. 
Then there exist $\lambda \in K$ and $Q \in \tilde{B}_n(R)$ such that 
$P = \lambda Q$.  
\end{lemma}

\prf Set $P = (p_{ij})$. 
Let $v$ be a valuation of $R$. 
We claim that $v(p_{11}) = v(p_{22}) = \cdots = v(p_{nn})$ 
and that $v(p_{ij}) \ge v(p_{11})$. 
From this claim, $\lambda:=p_{11}$ and 
$Q:=1/p_{11}\cdot P$ are what we want, and hence the statement can be proved. 

By the assumption, $P E_{ij} P^{-1} \in {\mathcal B}_n(R)$, and 
hence the $(i, j)$-entry $p_{ii}/p_{jj} \in R$. 
Since $P^{-1} {\mathcal B}_n(R) P 
={\mathcal B}_n(R)$ also holds, the $(i, j)$-entry $p_{jj}/p_{ii}$ 
of $P^{-1} E_{ij} P$ is contained in $R$. 
Therefore $v(p_{ii}) = v(p_{jj})$. 
For each $i<j$, the $(i, j)$-entry $p_{ij}/p_{jj}$ of $P E_{jj} P^{-1}$ 
is contained in $R$, which conclude that $v(p_{ij}) \ge v(p_{jj})=v(p_{11})$. 
Thereby we have proved the claim. 
\qed

\bigskip 

{\it Proof of Proposition \ref{prop-sep}.} 
We prove that ${\rm Ch}_n(\Gamma)_B$ is separated by using valuative criterion. 
Let $R$ be a valuation ring and $K$ be its quotient field. 
Suppose that $[ \rho ]$ and $[ \rho' ]$ be two $R$-valued points of 
${\rm Ch}_n(\Gamma)_B$ which coincide as $K$-valued points. 
We show that $[ \rho ] = [ \rho' ]$ as $R$-valued points. 
Let us take representatives $\rho, \rho' : \Gamma \to 
{\mathcal B}_n(R)$ of $[ \rho ], [ \rho' ]$, respectively.
Since $\rho, \rho'$ are equivalent to each other over $K$, there exists 
$P \in \tilde{B}_n(K)$ such that $P \rho P^{-1} = \rho'$. 
The algebra ${\mathcal B}_n(R)$ is generated by the image of 
$\rho$ or $\rho'$ over $R$, and hence $P {\mathcal B}_n(R) P^{-1} =
{\mathcal B}_n(R)$. 
By Lemma \ref{lemma-sep}, there exist $\lambda \in K$ and $Q \in \tilde{B}_n(R)$ 
such that $P=\lambda Q$. We obtain $Q \rho Q^{-1} = \rho'$, and we conclude that $\rho$ and $\rho'$ are equivalent over $R$.  
\qed

\bigskip

Summarizing the above discussion, we obtain: 

\begin{theorem}
Let $\Gamma$ be a group or a monoid. 
The sheafification with respect to Zariski topology of 
the following functor is representable by a separated scheme ${\rm Ch}_n(\Gamma)_B$ over 
${\Bbb Z}$: 
\[
\begin{array}{ccccl}
{\mathcal EqB}_n(\Gamma) & : & ({\bf Sch}) & \to & ({\bf Sets}) \\
 & & X & \mapsto & \{ \mbox{ rep. with Borel mold of deg $n$ for $\Gamma$ } \} /\sim 
\end{array}
\]
If $\Gamma$ is finitely generated, then ${\rm Ch}_n(\Gamma)_B$ 
is of finite type over ${\Bbb Z}$. 
\end{theorem}

\bigskip 

\begin{remark}\rm
In this paper, we deal with only representations of groups or monoids. 
However, we can construct the moduli of representations with Borel mold
for an arbitrary associative algebra. 
Let $A$ be an associative algebra over a commutative ring $R$. 
We define a representation with Borel mold for $A$ on a scheme 
over $R$ in a similar way as the group case. 
Then we can construct the moduli scheme of representations with Borel mold 
separated over $R$. If $A$ is a finitely generated algebra over $R$, then 
the moduli is of finite type over $R$.  
(The fact that the moduli is quasi-compact follows from that 
there exist a noetherian subring $S$ of $R$ and 
a finitely generated subalgebra $A_0$ of $A$ over $S$ such that 
$A_0 \otimes_{S} R \to A$ is surjective and 
the morphism ${\rm B}_n(A)_{B} \to {\rm B}_n(A_0)_{B}$ is affine 
and hence quasi-compact.)  
\end{remark}

%
%

\section{Basic results 
} 

In this section we introduce basic results on 
the moduli of representations with Borel mold. 

\bigskip

Let $\Gamma$ be a group or a monoid. 
Let $\rho$ be a representation with Borel mold for $\Gamma$ 
on a scheme $X$.
Let us define the action of $\Gamma$ on the trivial vector 
bundle ${\mathcal O}_X^{\oplus n}$ by 
$\Gamma \stackrel{\rho}{\to} {\rm M}_n(\Gamma(X, {\mathcal O}_X))
={\rm End}_{{\mathcal O}_X}({\mathcal O}_X^{\oplus n})$. 
Then we have the following proposition. 

\begin{proposition}\label{prop-flag}
For each $1 < i < n$, 
there exists a unique $\Gamma$-invariant subbundle 
${\mathcal E}_i \subseteq {\mathcal O}_X^{\oplus n}$ of 
rank $i$ on $X$. The $\Gamma$-invariant subbundles 
$0 \subseteq {\mathcal E}_1 \subseteq {\mathcal E}_2 
\subseteq \cdots \subseteq {\mathcal E}_{n-1} \subseteq 
{\mathcal O}_X^{\oplus n}$ form a complete flag of 
${\mathcal O}_X^{\oplus n}$. 
\end{proposition}

\prf
Lemma \ref{lemma-mold-sub} follows the 
uniqueness of $\Gamma$-invariant 
subbundles. 
Since we have $\Gamma$-invariant subbundles locally, 
by gluing them together 
we obtain unique $\Gamma$-invariant 
subbundles globally. 
\qed

\bigskip

From the above proposition, we easily obtain: 

\begin{theorem}\label{thpin}
The representation variety with Borel mold 
${\rm Rep}_n(\Gamma)_B$ has a 
unique complete flag of $\Gamma$-invariant 
subbundles $0 \subseteq {\mathcal E}_{\Gamma, n}^{(1)} 
\subseteq \cdots \subseteq 
{\mathcal E}_{\Gamma, n}^{(n-1)} \subseteq  
{\mathcal O}_{{\rm Rep}_n(\Gamma)_{B}}^{\oplus n}$ which has the 
universal property:
for any scheme $X$ and for any representation $\rho$ of degree $n$ 
with Borel mold for $\Gamma$
on $X$, the unique $\Gamma$-invariant subbundles  
${\mathcal E}_i$ of rank $i$ 
on $X$ is obtained as $f^{\ast}{\mathcal E}_{\Gamma, n}^{(i)}$,
where $f : X \to {\rm Rep}_n(\Gamma)_{B}$ is 
the morphism associated to $\rho$.
\end{theorem}

\bigskip

\begin{remark}\rm
For a representation with Borel mold $\rho$ for $\Gamma$ on 
a scheme $X$, we have unique $\Gamma$-invariant 
subbundles $0 \subseteq {\mathcal E}_1 \subseteq {\mathcal E}_2 
\subseteq \cdots \subseteq {\mathcal E}_{n-1} \subseteq 
{\mathcal O}_X^{\oplus n}$ on $X$ by Proposition
\ref{prop-flag}.
The action of $\Gamma$ on ${\mathcal E}_i/{\mathcal E}_{i-1}$ 
induces the character $\chi_i$ of $\Gamma$ for each $1 \le i \le n$. 
The correspondence $\rho \mapsto (\chi_1, \chi_2, \ldots, 
\chi_n)$ gives us a morphism ${\rm Ch}_n(\Gamma)_B \to 
{\rm Ch}_1(\Gamma)\times \cdots \times{\rm Ch}_1(\Gamma)$. 
Here ${\rm Ch}_1
(\Gamma)$ is the moduli of 
characters for $\Gamma$. 

In the case $n=2$, we have a morphism ${\rm Ch}_2(\Gamma)_B \to 
{\rm Ch}_1(\Gamma) \times {\rm Ch}_1(\Gamma)$. 
The fiber at $(\chi_1, \chi_2)$ is the projective space of the extension classes 
of characters $(\chi_1, \chi_2)$. 
However, in this article we will not go into details on the relation between the moduli 
of representations with Borel mold and the 
extension classes of characters.  

\end{remark}

\bigskip 

We denote by ${\rm Rep}_n(m)_B$ 
the representation variety with Borel mold 
for the free monoid of rank $m$. 
The following proposition follows that ${\rm Rep}_n(m)_B$ contains 
the representation variety with Borel mold 
${\rm Rep}_n({\rm F}_m)_B$ for 
the free group of rank $m$ as an open subscheme.  

\begin{proposition}
Let $\Upsilon_{m} = \langle \alpha_1, \ldots, \alpha_m \rangle$ 
be the free monoid of rank $m$. 
Let ${\rm F}_m = \langle \alpha_1, \ldots, \alpha_m \rangle$ 
be the free group of rank $m$.
The inclusion $\Upsilon_{m} \to {\rm F}_m$ by 
$\alpha_i \to \alpha_i$ induces an open immersion 
${\rm Rep}_n({\rm F}_m)_B \to {\rm Rep}_n(m)_B$.  
\end{proposition}

\prf
Restricting each representation $\rho$ with Borel mold for ${\rm F}_m$ 
to the free monoid $\Upsilon_m$, we can obtain a morphism 
${\rm Rep}_n({\rm F}_m)_B \to {\rm Rep}_n(m)_B$. 
Indeed, by the Cayley-Hamilton theorem, 
$\rho(\alpha_i^{-1})$ is expressed as a polynomial 
of $\rho(\alpha_i)$, and hence 
$\langle \rho(\alpha_1), \ldots, \rho(\alpha_m) \rangle$ 
generates a Borel mold. 
It is easy to check that 
the morphism ${\rm Rep}_n({\rm F}_m)_B \to {\rm Rep}_n(m)_B$ 
is an open immersion. 
\qed

\bigskip

In the case $n=1$, we see that ${\rm Rep}_1(m) \cong {\rm Rep}_1(m)_B \cong  {\mathbb A}_{\mathbb Z}^{m}$ and 
${\rm Rep}_1(F_m) \cong {\rm Rep}_1(F_m)_B \cong 
({\mathbb A}^{1}_{\mathbb Z} \setminus \{ 0 \} )^{m}$.  
In the case $n \ge 2$ and $m=1$, we also see that ${\rm Rep}_n(1)_B = \emptyset$ and 
${\rm Rep}_n(F_1)_B = \emptyset$.  
In the case $n \ge 2$ and $m \ge 2$, ${\rm Rep}_n(m)_B$ and 
${\rm Rep}_n({\rm F}_m)_B$ are non-empty (see \cite{topos}).
Furthermore we have: 

\begin{proposition}
Let $n \ge 2$ and $m \ge 2$. The scheme ${\rm Rep}_{n}(m)_{B}$ 
is smooth over ${\Bbb Z}$.
In particular, ${\rm Rep}_n({\rm F}_m)_B$ is smooth over
${\Bbb Z}$. 
\end{proposition}

\prf
Let $A$ be an artin local ring and let $I$ be 
an ideal of $A$ with $I^2 = 0$.
Let $\overline{\rho} \in {\rm Rep}_{n}({m})_{B}(A/I)$.
Then we show that there exists 
$\rho \in {\rm Rep}_{n}(m)_{B}(A)$ 
such that the reduction of $\rho$ to $A/I$ is 
equal to $\overline{\rho}$.
We take a system of free generators $\{ \alpha_{1}, \alpha_{2}, 
\ldots, \alpha_{m} \}$ of the free monoid $\Upsilon_{m}$.
There exists $\overline{P} \in 
{\rm GL}_{n}(A/I)$ such that 
$\overline{P}^{-1}\overline{\rho}(\alpha_{i})\overline{P}$
is an upper triangular matrix for each $1 \le i \le m$.
Take $P \in {\rm GL}_{n}(A)$ such that the reduction to $A/I$ is 
equal to $\overline{P}$.
We also take an upper triangular matrix $X_i 
\in {\rm GL}_n(A)$
as a lift of 
$\overline{P}^{-1}\overline{\rho}(\alpha_{i})\overline{P}$ 
for each $i$. 
Then we define the representation $\rho :
\Upsilon_{m} \to {\rm GL}_{n}(A)$ by 
$\rho(\alpha_{i}) := PX_{i} P^{-1}$ for each $1 \le i \le m$.
We easily see that $\rho$ is the desired representation with Borel mold, 
and hence we have proved the statement.
\qed

\bigskip 

\begin{corollary}
For $m \ge 2$, the moduli scheme of representations with Borel mold 
${\rm Ch}_n(m)_B$ 
for the free monoid $\Upsilon_m$ 
is smooth over ${\Bbb Z}$. 
In particular, the open subscheme ${\rm Ch}_n({\rm F}_m)_B$ of 
${\rm Ch}_n(m)_B$ is also smooth over ${\Bbb Z}$. 
\end{corollary}

\prf
The quotient morphism ${\rm Rep}_n(m)_B \to 
{\rm Ch}_n(m)_B$ gives a ${\rm PGL}_n$-principal 
fiber bundle. Since ${\rm Rep}_n(m)_B$ is smooth over 
${\Bbb Z}$, so is ${\rm Ch}_n(m)_B$.  
\qed

\begin{remark}\rm
In \cite{topos} we have proved that ${\rm Ch}_n(m)_B$ is 
a fibre bundle over the configuration space $F_n({\mathbb A}^m_{\mathbb Z})$ 
of the affine space ${\mathbb A}^m_{\mathbb Z}$ with fibre 
$({\mathbb P}^{m-2}_{\mathbb Z})^{n-1} \times ({\mathbb A}^{m-1}_{\mathbb Z})^{(n-2)(n-1)/2}$ 
with respect to Zariski topology. 
In particular, the rational function field of ${\rm Ch}_n(m)_B$ is rational over ${\mathbb Q}$ 
if $n=1$ or $n, m \ge 2$. Furthermore, if $k$ is a field, then 
${\rm Ch}_n(m)_B \otimes k$ is a smooth rational variety over $k$. 
\end{remark}

%
%

\section{The degree $2$ case}

In this section, we deal with representations of degree $2$ with Borel mold. 
In the degree $2$ case, a mold is a Borel mold if and only if 
it has rank $3$. 

\medskip 

The following proposition gives us one of characterizations 
of representations of degree $2$ with Borel mold 
for a group $\Gamma$.

\begin{proposition}\label{lemma-b-eq}
  Let $\Gamma$ be a group. 
Let $\rho$ be a representation of degree $2$ for $\Gamma$ on 
a scheme $X$. 
Then 
$\rho$ is a representation with Borel mold if and only if 
the following two conditions hold: 
\begin{enumerate}
\item\label{pin1} ${\rm tr}(\rho([\alpha, \beta ])) = 2$
for each $\alpha, \beta \in \Gamma$, 
where $[\alpha, \beta ] := \alpha\beta\alpha^{-1}\beta^{-1}$.
\item\label{pin2} the image of 
the composition $\Gamma 
\stackrel{\rho}{\to} {\rm GL}_2(\Gamma(X, {\mathcal O}_X)) 
\to {\rm GL}_2(k(x))$ 
is not an abelian group for each point $x \in X$, 
where $k(x)$ is the residue field of $x$.
\end{enumerate}
\end{proposition}

\prf
If $\rho$ is a representation with Borel mold, then 
for each $x \in X$ there exist a neighborhood $U$ of $x$ 
and $P \in {\rm GL}_2({\mathcal O}_X(U))$ such that 
$P^{-1}(\rho\!\mid_{U})(\gamma)P = 
\left(
  \begin{array}{cc}
    \ast & \ast \\
    0 & \ast 
  \end{array}
\right)$ 
for each $\gamma \in \Gamma$. 
Hence we have ${\rm tr}(\rho([\alpha, \beta])) = 2$ for each $\alpha, \beta 
\in \Gamma$. Since for the mold ${\mathcal O}_X[\rho(\Gamma)]
\otimes k(x)$ is a non-commutative algebra for $x \in X$, the condition
(\ref{pin2}) holds. 

Conversely suppose that 
two conditions (\ref{pin1}) and (\ref{pin2}) hold. 
Then we show that the subsheaf ${\mathcal O}_X[\rho(\Gamma)]$ 
of ${\rm M}_2({\mathcal O}_X)$ is a rank $3$ mold. 
If ${\mathcal O}_X[\rho(\Gamma)]$ is a rank $3$ mold, then 
it is a Borel mold from Corollary \ref{cormaxp}, which 
completes the proof. 
For each $x \in X$, there exists $\alpha, \beta \in \Gamma$ such that
$\rho(\alpha)$ and $\rho(\beta)$ are not commutative as 
elements of ${\rm GL}_2(k(x))$. 
From the assumption, the discriminant $\Delta(\rho(\alpha), \rho(\beta))$ 
in Definition \ref{def-appendix} is equal to $0$, 
since $\Delta(\rho(\alpha), \rho(\beta)) = 
\det(\rho(\alpha\beta))({\rm tr}(\rho([ \alpha, \beta])-2)$   
by Proposition \ref{prop-form-tr}. 
Proposition \ref{prop-appendix} follows that 
the subsheaf of ${\mathcal O}_U$-algebras 
${\mathcal O}_U[(\rho\!\mid_{U})(\alpha), (\rho\!\mid_{U})(\beta)]$ 
generated by $(\rho\!\mid_{U})(\alpha), (\rho\!\mid_{U})(\beta)$ 
is a rank $3$ mold 
on some affine neighbourhood $U$ of $x$.  
For each $\gamma \in \Gamma$, we only have to 
show that $(\rho\!\mid_{U})(\gamma) \in 
{\mathcal O}_U[(\rho\!\mid_{U})(\alpha), (\rho\!\mid_{U})(\beta)]$. 
By considering the subgroup of $\Gamma$ generated by $\alpha$, 
$\beta$, and $\gamma$, we have reduced to the case
that $\Gamma$ is a finitely generated 
group. 
Let $U = {\rm Spec}(R)$. 
Let us denote $\rho\!\mid_{U}$ by $\rho$. 
By changing $R$ to the subring of $R$ generated by 
all the entries of $(\rho\!\mid_{U})(\alpha), (\rho\!\mid_{U})(\beta)$, 
and $(\rho\!\mid_{U})(\gamma)$ over ${\Bbb Z}$ 
(if necessary, we may add more finitely many elements to the 
subring), 
we may assume that $R$ is a noetherian ring. 
Since we only need to prove that $(\rho\!\mid_{U})(\gamma) \in 
R[(\rho\!\mid_{U})(\alpha), (\rho\!\mid_{U})(\beta)]$ on a neighbourhood of 
each point of ${\rm Spec}R$, we may also 
assume that $R$ is a noetherian local ring. 


First suppose that $R$ is a reduced noetherian local ring. 
Since the subalgebra $R[\rho(\alpha), \rho(\beta)]$ of ${\rm M}_2(R)$ 
is a Borel mold, there exists $P \in {\rm GL}_2(R)$ such that 
$P^{-1}\rho(\alpha)P$ and $P^{-1}\rho(\beta)P$ generate 
the algebra ${\mathcal B}_2\otimes_{\Zah}R$. 
By changing $\rho$ to $P^{-1}\rho P$, we may assume that 
$\rho(\alpha)$ and $\rho(\beta)$ generate the algebra 
${\mathcal B}_2\otimes_{\Zah}R$. 
We can check that 
\begin{eqnarray}\label{eqn-b-ab}
\rho([\alpha, \beta]) 
=
\left(
\begin{array}{cc}
1 & u\\
0 & 1 
\end{array}
\right),
\end{eqnarray}
where $u \in R^{\times}$.
Put
\begin{eqnarray}\label{eqn-b-ab-2}
\rho(\gamma)  
 & = &
\left(
\begin{array}{cc}
a & b\\
c & d 
\end{array}
\right).
\end{eqnarray}
Then we have
\[
\rho([[\alpha, \beta], \gamma]) 
= 
\frac{1}{ad-bc}
\left(
\begin{array}{cc}
1 & u\\
0 & 1 
\end{array}
\right)
\left(
\begin{array}{cc}
a & b\\
c & d 
\end{array}
\right)
\left(
\begin{array}{cc}
1 & -u\\
0 & 1 
\end{array}
\right)
\left(
\begin{array}{cc}
d & -b\\
-c & a 
\end{array}
\right)
\]
\[
= \left(
\begin{array}{cc}
1 +\frac{acu + c^2u^2}{ad-bc}  & u - \frac{a^2u +acu^2}{ad-bc} \\
\frac{c^2u}{ad-bc} & 1 - \frac{acu}{ad-bc} 
\end{array}
\right).
\]
Since $\tr(\rho([[\alpha, \beta], \gamma])) = 2 + \frac{c^2u^2}{ad-bc}
= 2$, we have $c^2 = 0$. 
By the hypothesis that 
$R$ is reduced, we obtain $c = 0$, which 
implies that $\rho(\gamma) \in R[\rho(\alpha), \rho(\beta)]$.

Next we claim that 
if 
$(\rho(\gamma) \;{\rm mod}\; I) 
\in (R/I)[\rho(\alpha), \rho(\beta)]$ 
for an ideal $I$ of $R$ with $I^2 = 0$, then we can show that 
$\rho(\gamma) \in R[\rho(\alpha), \rho(\beta)]$. 
This claim and the result in the reduced case imply 
that $\rho(\gamma) \in R[\rho(\alpha), \rho(\beta)]$
for an arbitrary noetherian local ring $R$. 
As in the reduced case, we may assume that 
$\rho(\alpha)$ and $\rho(\beta)$ generate the algebra 
${\mathcal B}_2\otimes_{\Zah}R$ and that 
(\ref{eqn-b-ab}) and (\ref{eqn-b-ab-2}) hold. 
Since $(\rho(\gamma) \;{\rm mod}\; I)  
\in (R/I)[\rho(\alpha), \rho(\beta)]$, 
we have $c \in I$. 
By changing $\rho(\gamma)$ to $\rho(\gamma [\alpha, \beta ])$, 
if necessary,  
we may assume that the $(1, 2)$-entry $b$ of $\rho(\gamma)$ 
is contained in $R^{\times}$.  
Remark that at least one of $\rho(\alpha), \rho(\beta)$ 
is not commutative with $\rho(\gamma)$ as elements of 
${\rm GL}_2(R/m)$, where $m$ is the maximal ideal of $R$. 
Let $\langle \alpha, \beta \rangle$ be the subgroup of 
$\Gamma$ generated by $\alpha, \beta$. 

We see that there exists $\delta \in \langle \alpha, \beta \rangle$ 
such that 
$\rho(\delta)$ is not commutative with $\rho(\gamma)$ 
as elements of ${\rm GL}_2(R/m)$ and 
$\rho(\delta)$ 
has the form 
\[
\rho(\delta) = 
\left(
  \begin{array}{cc}
    p & q \\
    0 & r \\
  \end{array}
\right)
\]
with $p - r \in R^{\times}$. 
Indeed, suppose that there exists no such $\delta$. 
Then we have $\delta_1 \in \langle \alpha, \beta \rangle$
such that $\rho(\delta_1)$ is not commutative with $\rho(\gamma)$ 
as elements of ${\rm GL}_2(R/m)$ and 
$\rho(\delta_1)$ has the above form with $p - r \in m$. 
We also have $\delta_2 \in \langle \alpha, \beta \rangle$ 
such that $\rho(\delta_2)$ is commutative with $\rho(\gamma)$ 
as elements of ${\rm GL}_2(R/m)$ and 
$\rho(\delta_2)$ has the above form with $p - r \in R^{\times}$, 
because $\langle \rho(\alpha), \rho(\beta) \rangle$ generate 
${\mathcal B}_2 \otimes_{\mathbb Z} R/m$. 
Putting $\delta = \delta_1\delta_2$, we have such $\delta$. 

For such $\delta$, we obtain 
\[\displaystyle 
{\rm tr}(\rho([\gamma, \delta])) 
= 2 - c\frac{(p-r)\{ b(p-r) + q(d-a) \}}{(ad-bc)pr}. 
\] because $c \in I$ and $c^2 = 0$. 
Since $\rho(\gamma)$ and $\rho(\delta)$ are not commutative as 
elements of ${\rm GL}_2(R/m)$, 
$b(p-r) + q(d-a) \in R^{\times}$. 
From the assumption that ${\rm tr}(\rho([\gamma, \alpha])) = 2$, 
we have $c = 0$ because $p-r \in R^{\times}$. 
This implies that $\rho(\gamma) \in R[\rho(\alpha), \rho(\beta)]$. 
\qed

\bigskip 

From the above proposition, we have the following corollary. 

\begin{corollary}\label{cor-b-norm}
  Let $\rho$ be a representation of degree $2$ for a group $\Gamma$ on 
a scheme $X$. Suppose that 
$\rho$ satisfies the conditions (\ref{pin1}) and (\ref{pin2}) 
in Proposition \ref{lemma-b-eq}.
Then there exist an open covering $X = \cup_{i \in I}U_i$ and 
$P_i \in {\rm GL}_2(\Gamma(U_i, {\mathcal O}_X))$ 
such that 
\[
P_i^{-1}(\rho\!\mid_{U_i})(\gamma)P_i = 
\left(
  \begin{array}{cc}
    \ast & \ast \\
    0 & \ast 
  \end{array}
\right)
\]
for $\gamma \in \Gamma$ and 
$i \in I$. 
\end{corollary}


\begin{remark}\rm
The condition (\ref{pin1})  in Proposition 
\ref{lemma-b-eq} is necessary that 
$\rho$ can be normalized into 
upper triangular matrices as in Corollary \ref{cor-b-norm}.
Indeed, for such a representation $\rho$
we always have ${\rm tr}(\rho([\alpha, \beta])) =2$
for $\alpha, \beta \in \Gamma$.
Corollary \ref{cor-b-norm} does not always hold without the 
condition (\ref{pin2}) in Proposition \ref{lemma-b-eq}. 
The representation 
\[
\begin{array}{ccccl}
\rho & : & {\Bbb R} & \to & {\rm GL}_2({\Bbb R}) \\
 & & \theta & \mapsto & 
\left(
\begin{array}{cc}
\cos\theta & -\sin\theta \\
\sin\theta & \cos\theta
\end{array}
\right)
\end{array}
\]  
satisfies the condition (\ref{pin1}), but 
$\rho$ has no nontrivial invariant subspace 
in ${\Bbb R}^2$. Indeed, the above representation $\rho$ 
does not satisfy the condition (\ref{pin2}).
\end{remark}


\begin{remark}\rm
Note that any representation $\rho : \Gamma \to {\rm GL}_2({\Bbb F}_2)$
over the field ${\Bbb F}_2$ for a group $\Gamma$ is 
not a representation with Borel mold.
Indeed, if there exists a representation $\rho$ with Borel 
mold, then 
$\rho$ can be normalized into upper triangular matrices.
However the nonsingular upper triangular matrices over ${\Bbb F}_2$
are 
\[
\left(
\begin{array}{cc}
1 & 0 \\
0 & 1
\end{array}
\right) \mbox{ and }
\left(
\begin{array}{cc}
1 & 1 \\
0 & 1
\end{array}
\right),
\]
and hence $\rho(\Gamma)$ is an abelian group.
This is a contradiction.
\end{remark}


\bigskip 

From Theorem \ref{thpin}, we have a unique $\Gamma$-invariant 
sub-line bundle ${\mathcal L}_{\Gamma}$ of 
${\mathcal O}_{{\rm Rep}_{2}(\Gamma)_{B}}^{\oplus 2}$
on ${\rm Rep}_{2}(\Gamma)_{B}$.
We call ${\mathcal L}_{\Gamma}$ 
the {\it universal sub-line bundle} on 
${\rm Rep}_2(\Gamma)_B$.
Let us investigate the universal sub-line bundle 
${\mathcal L}_{\Gamma}$.

\begin{proposition}\label{propf2}
Let ${\rm F}_2 = \langle \alpha, \beta \rangle$ be the free group of rank $2$.
The universal sub-line bundle ${\mathcal L}_{{\rm F}_2}$
over ${\rm Rep}_{2}({\rm F}_2)_{B}$
is not trivial. 
\end{proposition}

\prf
We define the subvariety $X$ of ${\Bbb C}^4$
by $X := \{ (z_1, z_2, z_3, z_4) \in {\Bbb C}^4 \mid 
z_1^2+z_2^2+z_3^2+z_4^2 = 1 \}$.
The group ${\Bbb Z}/2{\Bbb Z} = \{ \pm 1 \}$ acts on 
${\Bbb C}^4$ and $X$ by 
$(z_1, z_2, z_3, z_4) \mapsto (-z_1, -z_2, -z_3, -z_4)$. 
We define $\varphi : X \to {\Bbb C}^2\setminus \{0\}$
by $(z_1, z_2, z_3, z_4) \mapsto (z_1 + \sqrt{-1}z_2, z_3 + \sqrt{-1}z_4)$.
We denote by $\psi$ the canonical projection 
${\Bbb C}^2 \setminus \{ 0 \} \to \Bbb{CP}^1$.
The morphisms $\varphi$ and $\psi$ induce the 
following diagram:
\[
\begin{array}{ccccc}
X & \stackrel{\varphi}{\to} & {\Bbb C}^2 \setminus \{ 0 \} & 
\stackrel{\psi}{\to} & \Bbb{CP}^1 \\
\downarrow & & \downarrow  & & \Vert \\
X/\{ \pm 1 \} & \stackrel{\overline{\varphi}}{\to} &
(\Bbb{C}^2 \setminus \{0\}) / \{ \pm 1 \} & \stackrel{\overline{\psi}}{\to} &
\Bbb{CP}^1.
\end{array}
\]
We denote by $f$ and $\overline{f}$ the compositions 
$\psi \circ \varphi$ and $\overline{\psi}\circ\overline{\varphi}$, 
respectively.
Let us consider the exact sequence 
\[ 0 \to {\mathcal O}_{\Bbb{CP}^1}(-1) \to 
{\mathcal O}_{\Bbb{CP}^1}^{\oplus 2} \to 
{\mathcal O}_{\Bbb{CP}^1}(1) \to 0\]
over $\Bbb{CP}^1$.
Taking the pull back by $\overline{f}$, we have the following exact sequence
\begin{equation}\label{exact1} 
0 \to\overline{f}^{\ast}{\mathcal O}_{\Bbb{CP}^1}(-1) \to 
{\mathcal O}_{X/\{ \pm 1 \}}^{\oplus 2} \to \overline{f}^{\ast}{\mathcal O}_{\Bbb{CP}^1}(1) \to 0.\end{equation}
We put ${\mathcal L} := \overline{f}^{\ast}{\mathcal O}_{\Bbb{CP}^1}(-1)$
and ${\mathcal M} := \overline{f}^{\ast}{\mathcal O}_{\Bbb{CP}^1}(1)$.
Note that ${\mathcal L} \cong {\mathcal M}^{-1}$.
We easily see that 
$\overline{\psi}^{\ast}{\mathcal O}_{\Bbb{CP}^1}(-1)^{\otimes 2}$
is trivial on $({\Bbb C}^{2}\setminus \{0\}) /\{ \pm 1 \}$, and hence 
we have ${\mathcal L}^{\otimes 2} \cong {\mathcal O}_{X/ \{ \pm 1 \}}$
and ${\mathcal M} \cong {\mathcal L}$.
The varieties $X$ and $X/\{ \pm 1 \}$ has the same homotopy types as 
$S^{3}$ and $\Bbb{RP}^3$, respectively.
We can regard the above diagram as follows up to homotopy:
\[
\begin{array}{ccc} 
S^3 & \stackrel{f}{\to}&  S^2 \\
\pi\downarrow & & \Vert \\
\Bbb{RP}^3 & \stackrel{\overline{f}}{\to} & S^2. 
\end{array}
\]
Here the map $\pi$ is the canonical projection. 
The map $f$ is the Hopf map.
The first Chern class $c_1({\mathcal O}_{{\Bbb C}{\Bbb P}^1}(-1))$ is a generator of 
${\rm H}^2(S^2, {\Bbb Z}) \cong {\Bbb Z}$
and $c_1({\mathcal  L}) \neq 0$
in ${\rm H}^2(\Bbb{RP}^3, {\Bbb Z}) \cong {\Bbb Z}/2{\Bbb Z}$.
Hence ${\mathcal L}$ is not trivial as a topological vector bundle.
Therefore we see that ${\mathcal L} \not\cong {\mathcal O}_{X/\{ \pm 1 \}}$.

Let us define a representation with Borel mold on $X/\{ \pm 1 \}$.
On the affine variety $X / \{ \pm 1 \}$, the exact sequence (\ref{exact1}) 
splits, and hence we have ${\mathcal O}_{X / \{ \pm 1 \}}^{\oplus 2} \cong 
{\mathcal L} \oplus {\mathcal M} \cong {\mathcal L} \oplus {\mathcal L}$.
We put $R := \Gamma(X / \{ \pm 1 \}, {\mathcal O}_{X / \{ \pm 1 \}})$.
Through the identification ${\rm M}_2(R) = 
{\rm End}_{{\mathcal O}_{X / \{ \pm 1 \}}}%
({\mathcal O}_{X / \{ \pm 1 \}}^{\oplus 2})
= {\rm End}_{{\mathcal O}_{X / \{ \pm 1 \}}} \\ 
({\mathcal L}\oplus {\mathcal M})$,
we define $\rho : {\rm F}_2 \to {\rm GL}_2(R)$ by 
\[
\rho(\alpha) :=
\left( 
\begin{array}{cc}
1 & 1\\
0 & 1
\end{array}
\right),
\rho(\beta) :=
\left( 
\begin{array}{cc}
\sqrt{-1} & 0\\
0 & -\sqrt{-1}
\end{array}
\right)\]
\[
\hspace*{20ex}
\in 
\left( 
\begin{array}{cc}
{\rm Hom}({\mathcal L}, {\mathcal L}) & {\rm Hom}({\mathcal M}, {\mathcal L})\\
{\rm Hom}({\mathcal L}, {\mathcal M}) & {\rm Hom}({\mathcal M}, {\mathcal M})
\end{array}
\right).
\]
Note that ${\rm Hom}({\mathcal L}, {\mathcal L}) =
{\rm Hom}({\mathcal L}, {\mathcal M}) = \cdots 
= {\mathcal O}_{X / \{ \pm 1 \}}$.
We see that $\rho$ is a representation with Borel mold and that 
${\mathcal L}$ is a unique $\Gamma$-invariant sub-line bundle.
For the morphism $g : X/ \{ \pm 1 \} \to {\rm Rep}_2({{\rm F}_2})_{B}$
associated to $\rho$, we have $g^{\ast}{\mathcal L}_{{\rm F}_2} 
= {\mathcal L} \not\cong {\mathcal O}_{X / \{ \pm 1 \}}$.
This implies  that ${\mathcal L}_{{\rm F}_2}$ is not trivial.
\qed

%

\bigskip

\begin{corollary}\label{cor-non-triv}
On the representation variety of 
degree $2$ with Borel mold ${\rm Rep}_2(2)_{B}$
for the free monoid of rank $2$, 
the universal sub-line bundle ${\mathcal L}_2$ is 
not trivial.  
\end{corollary}

\prf
Let $\Upsilon_2 = 
\langle \alpha, \beta \rangle$ be the free monoid of rank $2$.  
The inclusion $\Upsilon_2 \to
{\rm F}_2$ ($\alpha \mapsto \alpha, \;\beta \mapsto \beta$) 
induces the morphism $\phi: {\rm Rep}_2({\rm F}_2)_{B} 
\to {\rm Rep}_2(2)_{B}$. 
We easily see that $\phi^{\ast}{\mathcal L}_{2} =
{\mathcal L}_{{\rm F}_2}$. 
By the previous proposition the universal sub-line bundle 
${\mathcal L}_{{\rm F}_2}$ is not trivial, 
neither is ${\mathcal L}_{2}$. 
\qed

\bigskip


We show that the universal sub-line bundle ${\mathcal L}_{2}$ is
a $2$-torsion element of the Picard group.

\begin{proposition}\label{proptor}
For the universal sub-line bundle ${\mathcal L}_{2}$ 
on ${\rm Rep}_{2}(2)_{B}$, 
we have ${\mathcal L}_{2}^{\otimes 2} \cong {\mathcal O}_{{\rm Rep}_{2}(2)_{B}}$.
\end{proposition}

\prf
We denote ${\rm B}_2(\Upsilon_2)_B$ by ${\rm B}_2(2)_B$ for 
the free monoid $\Upsilon_2$. 
Let us consider the morphism 
\[
\begin{array}{ccccc}
\psi_{2}  & : & {\rm B}_{2}(2)_{B} \times {\rm GL}_2
& \to & {\rm Rep}_{2}(2)_{B} \\
& & (\rho, P)  & \mapsto & P^{-1}\rho P.
\end{array}
\]
We can easily see that 
$\psi_{2}$ is a smooth morphism. 
(Furthermore, $\psi_2$ is a principal fiber bundle with fiber 
$\tilde{B}_2 := 
\left\{
\left(
\begin{array}{cc}
\ast & \ast \\
0 & \ast
\end{array}
\right)
\right\} 
\subseteq {\rm GL}_2$
.)
We denote the universal nonsingular $2\times 2$ matrix on ${\rm GL}_2$ 
by $\left(
\begin{array}{cc}
p & q \\
r & s
\end{array}
\right)$.
We also denote the universal representations in 
${\rm B}_{2}(2)_{B}$ and 
${\rm Rep}_{2}(2)_{B}$ by
\[
\rho_{{\rm B}}(\gamma) =
\left(
\begin{array}{cc}
a(\gamma) & b(\gamma) \\
0 & d(\gamma)
\end{array}
\right)\]
and 
\[
\rho_{{\rm Rep}}(\gamma) =
\left(
\begin{array}{cc}
a'(\gamma) & b'(\gamma) \\
c'(\gamma) & d'(\gamma)
\end{array}
\right)\]
for $\gamma \in \Upsilon_2 = \langle \alpha, \beta \rangle$, respectively.
By putting $u:= (a(\alpha)-d(\alpha))b(\beta) - 
b(\alpha)(a(\beta)-d(\beta))$, we have 
\[
\rho_{{\rm B}}(\alpha)\rho_{{\rm B}}(\beta)
- \rho_{{\rm B}}(\beta)\rho_{{\rm B}}(\alpha)  = 
\left(
\begin{array}{cc}
0 & u \\
0 & 0
\end{array}
\right). 
\]
Since $\rho_{\rm B}$ is a representation with Borel mold, 
$u$ is an invertible global function on 
${\rm B}_{2}(2)_{B}$.
Hence ${\rm Ker}(\rho_{{\rm B}}(\alpha)\rho_{{\rm B}}(\beta)
- \rho_{{\rm B}}(\beta)\rho_{{\rm B}}(\alpha))$ 
is the universal sub-line bundle on 
${\rm B}_2(2)_{\rm B}$. 
Through the morphism $\psi_{2}$ we have 
\begin{multline*}
\rho_{{\rm Rep}}(\alpha)\rho_{{\rm Rep}}(\beta)
- \rho_{{\rm Rep}}(\beta)\rho_{{\rm Rep}}(\alpha)  \\ 
= 
\left(
\begin{array}{cc}
p & q \\
r & s
\end{array}
\right)^{-1}
(\rho_{\rm B}(\alpha)\rho_{\rm B}(\beta) 
-\rho_{\rm B}(\beta)\rho_{\rm B}(\alpha))
\left(
\begin{array}{cc}
p & q \\
r & s
\end{array}
\right) \\
 = 
\frac{1}{\Delta} \left(
\begin{array}{cc}
rsu & s^2u  \\
-r^2u & -rsu 
\end{array}
\right), \hspace*{1cm}
\end{multline*}
where $\Delta := ps-qr$.
Since ${\rm Ker} (\rho_{{\rm Rep}}(\alpha)\rho_{{\rm Rep}}(\beta)
- \rho_{{\rm Rep}}(\beta)\rho_{{\rm Rep}}(\alpha))$ 
is equal to the universal sub-line bundle ${\mathcal L}_{2} \subseteq 
{\mathcal O}_{{\rm Rep}_{2}(2)_{B}}^{\oplus 2}$, 
${}^{t}(-s^2u/\Delta, rsu/\Delta)$ and ${}^{t}(rsu/\Delta, -r^2u/\Delta)$ 
are global sections of ${\mathcal L}_{2}$. 
We define the prime divisors $D_1$ and $D_2$ on 
${\rm Rep}_{2}(2)_{B}$ by 
$D_1 := \psi_{2}(\{ r =0 \})$ and 
$D_2 := \psi_{2}(\{ s =0 \})$, 
respectively.
Let us denote the $(1, 2)$-entry and $(2, 1)$-entry of 
$\rho_{{\rm Rep}}(\alpha)\rho_{{\rm Rep}}(\beta)
- \rho_{{\rm Rep}}(\beta)\rho_{{\rm Rep}}(\alpha)$ 
by $b'$ and $c'$, respectively.  
Then $D_1 = \{ c' = 0\}$
and $D_2 = \{ b' = 0\}$ set-theoretically.
From the two global sections above  we see that 
${\mathcal L}_{2} \sim D_1 \sim D_2$.
Because $2\cdot D_1 \sim {\rm div}(c') \sim 0$,
we conclude that 
${\mathcal L}_{{2}}^{\otimes 2} \cong 
{\mathcal O}_{{\rm Rep}_{2}({2})_{B}}$. 
\qed 

\bigskip

\begin{corollary}\label{cortor}
Let ${\mathcal L}_{{\rm F}_2}$ be the universal sub-line bundle 
on ${\rm Rep}_2({\rm F}_2)_{\rm B}$. 
Then ${\mathcal L}_{{\rm F}_2}^{\otimes 2} \cong 
{\mathcal O}_{{\rm Rep}_2({\rm F}_2)_{\rm B}}$.  
\end{corollary}

\prf
The statement follows from that the pull-back of 
${\mathcal L}_2$ by the morphism ${\rm Rep}_2({\rm F}_2)_{\rm B} 
\to {\rm Rep}_2(2)_{\rm B}$ is equal to ${\mathcal L}_{{\rm F}_2}$. 
\qed

\bigskip

From Corollary \ref{cortor}
we have the following: 

\begin{corollary}\label{cor-grt}
Let $\Gamma$ be a group generated by two elements.
Let $R$ be a commutative ring such that  
${\rm Pic}({\rm Spec}(R))$ has no $2$-torsion element.
For each representation $\rho : \Gamma \to
{\rm GL}_2(R)$ with Borel mold, we have some $P \in {\rm GL}_2(R)$
such that \[P\rho(\gamma) P^{-1} =
\left(
\begin{array}{cc}
a(\gamma) & b(\gamma) \\
0 &d(\gamma)
\end{array}
\right) \] 
for each $\gamma \in \Gamma$.
\end{corollary}

\prf
Let us consider a closed immersion 
$f: {\rm Rep}_2(\Gamma)_{\rm B} \to {\rm Rep}_2({\rm F}_2)_{\rm B}$ 
induced by a surjective morphism ${\rm F}_2 \to \Gamma$.  
From Corollary \ref{cortor}, we see that   
$f^{\ast}{\mathcal L}_{{\rm F}_2} = {\mathcal L}_{\Gamma}$ is 
a $2$-torsion element of the Picard group. 
Hence the pull-back of ${\mathcal L}_{\Gamma}$ on $R$ is trivial, 
which follows the claim. 
\qed 

\bigskip

Let us discuss the free group of rank $\ge 3$ case.
 
\begin{proposition}\label{propf3}
For the free group ${\rm F}_{m}$ with $m \ge 3$, 
we have ${\mathcal L}_{{\rm F}_m}^{\otimes n} \not\cong
{\mathcal O}_{{\rm Rep}_{2}({\rm F}_{m})_{{B}}}$ 
for each integer $n \neq 0$. 
\end{proposition}

\prf
There exists an affine scheme $X$ (over ${\Bbb C}$) which satisfies 
the following condition:
$X$ has a sub-line bundle 
${\mathcal L} \subseteq {\mathcal O}_{X}^{\oplus 2}$
such that ${\mathcal L}^{\otimes n} \not\cong {\mathcal O}_{X}$
for each integer $n \neq 0$.
Example \ref{exsch} gives us such an affine scheme $X$.
We denote ${\mathcal O}_{X}^{\oplus 2}/{\mathcal L}$ 
by ${\mathcal M}$. Then we have ${\mathcal O}_{X}^{\oplus 2} \cong
{\mathcal L} \oplus {\mathcal M}$ and 
${\mathcal M} \cong {\mathcal L}^{-1}$.
Since ${\mathcal L}$ is generated by two global sections,
${\rm Hom}({\mathcal M}, {\mathcal L}) \cong {\mathcal L}^{\otimes 2}$
is generated by some two global sections.
Suppose that $s, t \in {\rm Hom}({\mathcal M}, {\mathcal L})$
are global sections which generates ${\rm Hom}({\mathcal M}, {\mathcal L})$.
For ${\rm F}_3 = \langle \alpha, \beta, \gamma \rangle$, we define 
a representation $\rho_3$ on $X$ by
\[
\rho_3(\alpha):=\left(
\begin{array}{cc}
\sqrt{-1} & 0 \\
0 & -\sqrt{-1}
\end{array}
\right),
\rho_3(\beta):=\left(
\begin{array}{cc}
1 & s \\ 
0 & 1 
\end{array} 
\right), \rho_3(\gamma):=\left(
\begin{array}{cc}
1 & t \\
0 & 1
\end{array}
\right)
\]
\[
\in 
{\rm End}_{{\mathcal O}_{X}}({\mathcal O}_{X}^{\oplus 2}) =
\left(
\begin{array}{cc}
{\rm Hom}({\mathcal L}, {\mathcal L}) & 
{\rm Hom}({\mathcal M}, {\mathcal L}) \\
{\rm Hom}({\mathcal L}, {\mathcal M}) & 
{\rm Hom}({\mathcal M}, {\mathcal M})
\end{array}
\right).
\]
For the free group ${\rm F}_m$ with $m \ge 4$,
we define a representation $\rho_{m}$ on $X$
by taking the composite of a surjection ${\rm F}_m 
\twoheadrightarrow {\rm F}_3$
and the above $\rho_{3}$.
We can check that $\rho_{m}$ is a representation with Borel mold.
For the morphism $g : X \to {\rm Rep}_2({\rm F}_m)_{B}$ 
associated to $\rho_m$, we have  
$g^{\ast}{\mathcal L}_{{\rm F}_m} \cong {\mathcal L}$.
This implies that ${\mathcal L}_{{\rm F}_m}^{\otimes n}
\not\cong {\mathcal O}_{{\rm Rep}_2({\rm F}_m)_{B}}$
for each integer $n \neq 0$.
\qedd 

\bigskip

\begin{corollary}
For the universal sub-line bundle ${\mathcal L}_m$ for 
the free monoid $\Upsilon_{m}$ with $m \ge 3$, 
we have ${\mathcal L}_{m}^{\otimes n} \not\cong
{\mathcal O}_{{\rm Rep}_{2}(m)_{{B}}}$ 
for each integer $n \neq 0$. 
\end{corollary}

\prf
We can prove the statement in the same way as Corollary \ref{cor-non-triv}. 
\qed

\bigskip

The following example has been used in the proof of 
the previous proposition. 

\begin{example}\rm\label{exsch}
Let us consider the affine variety $X := \{ (z_1, z_2, z_3) \in {\Bbb C}^{3}
\mid z_1^2 + z_2^2 + z_3^2 = 1 \}$.
The variety $X$ has the same homotopy type as $S^2$.
We define the morphism $f : X \to \Bbb{CP}^1$
by $(z_1, z_2, z_3) \mapsto (z_1 + \sqrt{-1}z_2, z_3)$.
Then the morphism $f$ can be regarded as 
a degree $2$ map $S^2 \to S^2$ up to homotopy.
We define the sub-line bundle ${\mathcal L}$ on 
$X$ by ${\mathcal L} := f^{\ast}{\mathcal O}_{\Bbb{CP}^1}(-1) \subset {\mathcal O}_X^{\oplus 2}$.
Since $c_1({\mathcal L}) \neq 0 \in {\rm H}^{2}(X, {\Bbb Z}) \cong {\Bbb Z}$,
the line bundle ${\mathcal L}$ is not a torsion 
in the Picard group ${\rm Pic}(X)$.
\end{example}

\begin{remark}\rm
From Propositions \ref{propf2} and \ref{proptor}
we see that the $2$-torsion part of the cohomology 
does not vanish:
${\rm H}^{2}%
({\rm Rep}_{2}({\rm F}_{2})_{B}\otimes_{{\Bbb Z}} {\Bbb C}, {\Bbb Z})_{2}
\neq 0$.  
By Proposition \ref{propf3} we also see that 
the free part of ${\rm H}^{2}%
({\rm Rep}_{2}({\rm F}_{m})_{B}\otimes_{{\Bbb Z}} {\Bbb C}, {\Bbb Z})$
does not vanish
for each $m \ge 3$.
The topology of the representation varieties and the 
character varieties with Borel mold has been investigated in 
\cite{topos}, \cite{rational}, and \cite{notestopos}.
\end{remark}

\section{Appendix}

In this appendix, we prove several propositions 
on discriminants in the degree $2$ case. 
These propositions have been used in the previous section. 
Discriminants are invariants which describe 
open subschemes of absolutely irreducible 
representations in the representation varieties. More 
precisely, see \cite{kj-Saito93}, \cite{Nkmt00}, and \cite{Nkmt-deg2}. 

\begin{definition}[\cite{kj-Saito93}]\label{def-appendix}\rm
Let $R$ be a commutative ring.

For $A, B \in \mat{2}{R}$ we define the {\it discriminant} 
$\Delta(A, B)$ by 
\[\begin{array}{ccl}
\Delta(A, B) & := & \tr(A)^2\det(B) + \tr(B)^2\det(A) 
+ \tr(AB)^2 \\ 
& & \hspace{10ex} - \tr(A)\tr(B)\tr(AB) -4\det(A)\det(B). 
\end{array}
\]
From the definition we see that 
$\Delta(A, B) = \Delta(B, A)$.
\end{definition}

\begin{remark}\label{remarktwodeltas}\rm
The discriminant $\Delta(A, B)$ above 
is closely related to the discriminant in 
\cite{Nkmt00}.
For $A, B, C, D \in \mat{2}{R}$ we define the {\it discriminant} 
of degree $2$ in \cite{Nkmt00}
by 
\[
\Delta(A, B, C, D) := \det \left(
\begin{array}{cccc}
a_1 & a_2 & a_3 & a_4 \\
b_1 & b_2 & b_3 & b_4 \\
c_1 & c_2 & c_3 & c_4 \\
d_1 & d_2 & d_3 & d_4 
\end{array}
\right),
\]
where 
\[
A=\left(
\begin{array}{cc}
a_1 & a_2 \\
a_3 & a_4
\end{array}
\right),
B=\left(
\begin{array}{cc}
b_1 & b_2 \\
b_3 & b_4
\end{array}
\right),
C=\left(
\begin{array}{cc}
c_1 & c_2 \\
c_3 & c_4
\end{array}
\right),
D=\left(
\begin{array}{cc}
d_1 & d_2 \\
d_3 & d_4
\end{array}
\right).
\]
Note that $\Delta(A, B, C, D) \in R^{\times}$ 
if and only if $\{ A, B, C, D \}$ is an $R$-basis of $\mat{2}{R}$.
For $A, B \in {\rm M}_2(R)$, we have 
$\Delta(A, B)=-\Delta(I_2, A, B, AB)$.
\end{remark}

\bigskip 

We can easily obtain 
the following proposition.

\begin{proposition}[cf. \cite{kj-Saito93}]\label{prop-form-tr}
For each $A, B \in {\rm GL}_2(R)$, we have 
\[
\Delta(A, B)= \det(AB)({\rm tr}(ABA^{-1}B^{-1})-2). 
\]
\end{proposition}

\bigskip 

The following proposition has been used in Proposition \ref{lemma-b-eq}. 

\begin{proposition}\label{prop-appendix}
Let $A, B \in {\rm M}_2(R)$. 
Suppose that 
$AB \neq BA$ as matrices in ${\rm M}_2(k(\wp))$ for each prime ideal $\wp \in {\rm Spec} R$, 
where $k(\wp) := R_{\wp}/\wp R_{\wp}$.  
If $\Delta(A, B) = 0$, then   
the $R$-subalgebra $R[ A, B ]$ of ${\rm M}_2(R)$ generated by 
$A$ and $B$ is contained in the $R$-submodule $R\cdot I_2+ R\cdot A+R\cdot B$. 
\end{proposition}

\prf
First we show the claim that $AB$ is expressed as a linear combination of 
$\{ I_2, A, B \}$. 
For proving this, we can assume that 
$R$ is a local ring. 
Indeed, let us define the ideal $J$ of $R$ by 
\[J:= \{ a \in R \mid a AB \mbox{ is expressed as a 
linear combination of $I_2, A, B$ } \}.  
\] 
If the claim is true for the local ring case, then  
$AB$ is expressed as a linear combination of 
$I_2, A, B$ in ${\rm M}_2(R_{\wp})$ 
for each prime ideal $\wp \in {\rm Spec}R$.
Hence $J \not\subset \wp$ for each $\wp \in {\rm Spec}R$, 
which implies $J=R$. 
Since $1 \in J$, the claim is true for an arbitrary ring $R$. 

Assume that $(R, m, k)$ is a local ring. 
From the hypothesis, 
$AB \neq BA$ as matrices in ${\rm M}_2(k)$.  
Then $I_2, A$, and $B$ are linearly independent in ${\rm M}_2(k)$, and hence 
they are also linearly independent in ${\rm M}_2(R)$.   
By Remark \ref{remarktwodeltas}, we have   
\[
\Delta(A, B) = -\Delta(I_2, A, B, AB) = -\det \left(
\begin{array}{cccc}
1 & 0 & 0 & 1 \\ 
a_1 & a_2 & a_3 & a_4 \\
b_1 & b_2 & b_3 & b_4 \\
c_1 & c_2 & c_3 & c_4 \\
\end{array}
\right) = 0, 
\]
where 
\[
A=\left(
\begin{array}{cc}
a_1 & a_2 \\
a_3 & a_4
\end{array}
\right),
B=\left(
\begin{array}{cc}
b_1 & b_2 \\
b_3 & b_4
\end{array}
\right),
AB=\left(
\begin{array}{cc}
c_1 & c_2 \\
c_3 & c_4
\end{array}
\right). 
\] 
At least one of the $3\times 3$ minor determinants  
of $T := \left(
\begin{array}{cccc}
1 & 0 & 0 & 1 \\ 
a_1 & a_2 & a_3 & a_4 \\
b_1 & b_2 & b_3 & b_4 \\
\end{array}
\right)$ is a unit of $R$. 
If $P := \left(
\begin{array}{ccc}
1 & 0 & 0  \\ 
a_1 & a_2 & a_3  \\
b_1 & b_2 & b_3  \\
\end{array}
\right) \in {\rm GL}_3(R)$, for example, then 
\begin{eqnarray}\label{eq:lastmatrix}  
\left( 
\begin{array}{cc} 
P^{-1} & 0 \\
0        & 1 \\
\end{array}
\right) 
\left(
\begin{array}{cccc}
1 & 0 & 0 & 1 \\ 
a_1 & a_2 & a_3 & a_4 \\
b_1 & b_2 & b_3 & b_4 \\
c_1 & c_2 & c_3 & c_4 \\
\end{array}
\right) = 
\left(
\begin{array}{cccc}
1 & 0 & 0 & \ast \\ 
0 & 1 & 0 & \ast \\
0 & 0 & 1 & \ast \\
c_1 & c_2 & c_3 & c_4 \\
\end{array}
\right).  
\end{eqnarray}  
Let $v(X) := (x_1, x_2, x_3, x_4) \in R^4$ for 
$X = \left( 
\begin{array}{cc}
x_1 & x_2 \\
x_3 & x_4 \\
\end{array} 
\right) \in {\rm M}_2(R)$. 
The right hand of (\ref{eq:lastmatrix}) has 
the form ${}^{t}\left( v(X_1), v(X_2), v(X_3), v(AB) \right) 
\in {\rm M}_4(R)$ for some  
$X_1, X_2, X_3 \in R\cdot I_2 + R\cdot A + R\cdot B$.  
The determinant of the matrix (\ref{eq:lastmatrix}) is 0.  
By matrix row transformation, we see that 
$v(AB)$ is expressed as a linear combination of 
$\{ v(X_1), v(X_2), v(X_3) \}$ and hence that 
$AB$ is expressed as a linear combination of 
$\{ I_2, A, B \}$. In the case that another $3\times 3$ minor determinant 
of $T$ is a unit of $R$, we also see that $AB \in R\cdot I_2 + R\cdot A + R\cdot B$. 
Similarly we can prove that $BA$ is expressed as a linear combination of 
$\{ I_2, A, B \}$ because $\Delta(A, B)=\Delta(B, A)=-\Delta(I_2, A, B, BA) = 0$. 

Next we show that any monomial of $A$ and $B$ is contained in the 
$R$-submodule $R\cdot I_2 + R\cdot A + R\cdot B$. 
Each monomial of length $\ge 2$ contains one of 
$AA, AB, BA$, and $BB$ as a subsequence. 
By the Cayley-Hamilton theorem and the above discussion, 
we can reduce $AA, AB, BA$, and $BB$ 
to monomials of length one.  
Hence each monomial is contained in 
$R\cdot I_2 + R\cdot A + R\cdot B$ by induction. 
Thus we have completed the proof. 
\qed

\begin{remark}\rm 
Proposition~\ref{prop-appendix} does not hold 
without the hypothesis that $AB \neq BA$ in ${\rm M}_2(k(\wp))$. 
Indeed, let $R := {\Bbb C}[\varepsilon]/(\varepsilon^4)$. Set 
$A = \varepsilon \left( 
\begin{array}{cc} 
0 & 1 \\
0 & 0 \\
\end{array} 
\right)$ and 
$B = \varepsilon \left( 
\begin{array}{cc} 
0 & 0 \\
1 & 0 \\
\end{array} 
\right)$. Then $AB = \varepsilon^2 \left( 
\begin{array}{cc} 
1 & 0 \\
0 & 0 \\
\end{array} 
\right)$ can not be expressed as a linear combination of 
$\{ I_2, A, B \}$ in ${\rm M}_2(R)$, although $\Delta(A, B) = 0$. 
This is a counterexample. 
\end{remark}

\bibliographystyle{amsplain}

\end{document}